\documentclass[12pt]{article}
\usepackage{amscd}
\usepackage{mathrsfs}
\usepackage{latexsym,amsfonts}
\usepackage{graphicx}
\usepackage[colorlinks,linkcolor=red,citecolor=blue]{hyperref}
\usepackage{amsmath,amscd,stmaryrd,amssymb,array, amsfonts,amsmath,amssymb}
\usepackage{verbatim}
\usepackage{enumitem}

\setenumerate[1]{itemsep=0pt,partopsep=0pt,parsep=\parskip,topsep=5pt}
\setitemize[1]{itemsep=0pt,partopsep=0pt,parsep=\parskip,topsep=0pt}
\setdescription{itemsep=0pt,partopsep=0pt,parsep=\parskip,topsep=5pt}
\usepackage{tocbibind}
\usepackage{xcolor}

\numberwithin{figure}{section}
 \numberwithin{equation}{section}

\newtheorem{theorem}{Theorem}[section]
\newtheorem{proposition}[theorem]{Proposition}
\newtheorem{definition}[theorem]{Definition}

\newtheorem{lemma}[theorem]{Lemma}

\newcommand{\abs}[1]{\lvert#1\rvert}
\newcommand{\dif}{\,\mathrm{d}}
\DeclareMathOperator{\supp}{supp}
\DeclareMathOperator{\diag}{diag}
\DeclareMathOperator{\dist}{dist}
\DeclareMathOperator{\Tr}{Tr}
\DeclareMathOperator{\Prj}{Prj}
\DeclareMathOperator{\Div}{div}

\newcommand{\bbR}{{\mathbb R}}
\newcommand{\cN}{{\mathcal N}}
 \newcommand{\cT}{{\mathcal T}}

\newcommand{\cS}{{\mathcal S}}

\usepackage{lipsum}
\newcommand\blfootnote[1]{%
\begingroup
\renewcommand\thefootnote{}\footnote{#1}%
\addtocounter{footnote}{-1}%
\endgroup
}

\def\Vs{\vskip8pt}\def\vs{\vskip4pt}

\begin{document}

\begin{center}
{\bf\large Saddle solutions for the fractional Choquard equation}
\end{center}

\vs\centerline{Ying-Xin Cui 
}  
\begin{center}
{\footnotesize
{Center for Applied Mathematics,  Tianjin University,\\
          Tianjin 300072,  China\\

{\em E-mail}:  cuiyingxin1993@163.com}}
\end{center}

\vs\centerline{Jiankang Xia$^\dag$\blfootnote{$\dag$Corresponding author.}
}  
\begin{center}
{\footnotesize
{School of Mathematics and Statistics, Northwestern Polytechnical University,\\
        Xi'an 710129, China\\

{\em E-mail}:  jiankangxia@nwpu.edu.cn}}
\end{center}
\Vs

\begin{minipage}{12cm}

{\footnotesize
{\bf Abstract.}We study the saddle solutions for the fractional Choquard equation
 \begin{align*}
(-\Delta)^{s}u+ u=(K_{\alpha}\ast|u|^{p})|u|^{p-2}u, \quad x\in \mathbb{R}^N
 \end{align*}
 where $s\in(0,1)$, $N\geq 3$ and $K_\alpha$ is the Riesz potential with order $\alpha\in (0,N)$.
For every Coxeter group $G$ with rank $1\leq k\leq N$ and $p\in[2,\frac{N+\alpha}{N-2s})$, we construct a $G$-saddle solution with prescribed symmetric nodal configurations. This is a counterpart for the fractional Choquard equation of saddle solutions to the Choquard equation and further completes the existence of non-radial sign-changing  solutions for this doubly nonlocal equation.

\medskip

}
\end{minipage}


\section{Introduction and main results}

\quad\, We consider the  Choquard equation involving a fractional Laplacian:
\begin{equation}\label{eqfrc}
\begin{aligned}
(-\Delta)^{s}u+ u=(K_{\alpha}\ast\abs{u}^{p})\abs{u}^{p-2}u, \quad x\in \mathbb{R}^N,
\end{aligned}\end{equation}
where $N\geq 3$, $K_{\alpha}:\mathbb R^N\to\mathbb R$ with $\alpha\in(0,N)$ is the Riesz potential, defined by
\begin{displaymath}\begin{aligned}
K_{\alpha}(x)=\frac{A_{\alpha}}{\abs{x}^{N-\alpha}}, \;\text{ with }\; A_{\alpha}=\frac{\Gamma(\frac{N-\alpha}{2})}{\Gamma(\frac{\alpha}{2})\pi^{\frac{N}{2}}2^{\alpha}}.
\end{aligned}\end{displaymath}
Here the fractional Laplacian $(-\Delta)^{s}$ with $s\in(0,1)$ is defined by
$$\mathcal{F}((-\Delta )^s u)(\xi)=\abs{\xi}^{2s}\mathcal{F}(u)(\xi),$$
where $\mathcal F$ denotes the Fourier transform. When $u$ is sufficiently smooth, the fractional Laplacian can also be expressed by
\begin{displaymath}
\begin{aligned}
 (-\Delta)^{s} u(x)=C_{N,s}\lim\limits_{\epsilon\rightarrow 0}\int_{\mathbb{R}^N \backslash B_{\varepsilon}(x)}
 \frac{u(x)-u(y)}{\abs{x-y}^{N+2s}}\dif y,
  \end{aligned}\end{displaymath}
with $C_{N,s}>0$ being a normalization constant. In the last few years, the fractional Laplacian  arises in the description of anomalous diffusion \cite{MK2000} and is treated as the infinitesimal generators of L\'evy stable diffusion processes \cite{L}, and since then various fractional equations were derived in distinct fields: game theory \cite{Caf}, minimal surfaces \cite{CSS} and finance \cite{CT}, to name a few.
 In the remarkable seminal work \cite{CS} of Caffarelli and Silvestre, the $s$-harmonic extension technique was introduced which makes it possible to transform the nonlocal problem into a local one via the Dirichlet--Neumann map.  When $s=1$, equation  \eqref{eqfrc}  is reduced to the Choquard equation or the Choquard--Pekar equation,
\begin{equation}\label{1.5}
\begin{aligned}
-\Delta u+ u=(K_{\alpha}\ast\abs{u}^{p})\abs{u}^{p-2}u, \quad x\in \mathbb{R}^N.
\end{aligned}\end{equation}
This type of nonlocal equation was adopted in many physical contexts such as the Hartree--Fock theory of one-component plasma \cite{Lieb}, the quantum mechanics of a polaron at rest \cite{P}, the self-gravitating matter model \cite{Pe} and so on.

The study of the Choquard equation \eqref{1.5} has been going on for many years. The existence and uniqueness of the groundstate solution was first considered by Lieb in 1976 \cite{Lieb} for the Choquard--Pekar equation \eqref{1.5}  with $N=3$, $p=2$ and $\alpha=2$. Lions \cite{Li} later showed the existence of infinitely many radially solutions to  the same model. Ma and Zhao \cite{MZ} proved that every positive solution of \eqref{1.5} is radially symmetric and monotone decreasing about some point under some assumptions on $N$, $\alpha$ and $p$. These restrictions used in \cite{MZ} were ultimately eliminated by Moroz and Van Schaftingen in \cite{MoS} and some qualitative properties of the groundstate were also investigated including positivity, radial symmetry and decay behaviors. Since the Choquard equation \eqref{1.5} has been widely investigated by employing variational methods, we limit ourselves to citing a few references \cite{MoS2,MoS4,SX1} and refer to their bibliographies for a broader list.

As for the fractional Choquard equation \eqref{eqfrc}, it turns to be a doubly nonlocal equation due to the appearance of the convolution and the fractional Laplacian. When $s=\frac{1}{2}$, Frank and Lenzmann \cite{Le2009} proved the radial symmetry, uniqueness and positivity of the $L^2$-critical groundstate solution in
the model of dynamics of pseudo-relativistic boson stars.  For the general setting,  d'Avenia, Siciliano and Squassina \cite{DSS} showed the existence of positive groundstate solutions under the assumptions of $p$ satisfying
\begin{equation}\label{1.3}\begin{aligned}
\frac{N+\alpha}{N}<p<\frac{N+\alpha}{N-2s}.
\end{aligned}\end{equation}
In addition, the radial symmetry, Morse index and asymptotic decays for the groundstate
solutions were investigated as well, see also \cite{BBMP}. They also showed that the above range of $p$ turns to be optimal for the existence of finite energy solutions so that the endpoints of the above interval are critical exponent with respect to the Hardy--Littlewood--Sobolev inequality
 \cite{LL}.  Shen, Gao and Yang \cite{SGY} obtained the groundstates for the fractional Choquard equation \eqref{eqfrc} with general nonlinearities in spirit of Berestycki--Lions. Some other results were also established by using variational methods on the existence and multiple solutions to the fractional Choquard equation \eqref{eqfrc}, we refer to \cite{CT2019, CL,GuoHu,LKZ,Luo,MZ2017} and the references therein.

In recent years, the study of the sign-changing solutions to the Choquard equation \eqref{1.5} has received a lot of attention. As we all know, the compactness embedding plays an important role in searching the entire nodal solutions for a class of variational nonlinear equations. As a result, the existence of sign-changing solutions of equation \eqref{1.5} has been established in
\cite{GG,HYY,SX1} under a compact setting, and we refer to \cite{ZW} for the existence of infinitely many radially sign-changing solutions to the fractional Choquard equation \eqref{eqfrc}. Recently,
 Ghimenti and Van Schaftingen  constructed in \cite{GS} a surprising odd solution with exactly two half-space nodal domains. Such an odd solution demonstrates different features of Choquard equation because its nodal set consists of  hyperplane whereas it is impossible for the nonlinear Schr\"odinger equations $-\Delta u+u=\abs{u}^{p-2}u$. As a consequence, by employing a minimax procedure on the Nehari nodal set, Ghimenti and Van Schaftingen \cite{GS} showed the existence of the minimal nodal solution of the Choquard equation \eqref{1.5} for $p\in(2, \frac{N+\alpha}{N-2})$, see
 \cite{GMS} for the case of $p=2$. This result is quite surprising since the energy of the minimal nodal solution is strict less than twice the groundstate energy, while the nonlinear Schr\"odinger equation has the doubling energy property, see e.g., \cite{Weth2006}. Moreover, it turns out some degeneracy for the case $p<2$ which leads to the minimizer on the Nehari nodal set is still the groundstate solution.
Wang and Xia \cite{XW2,XW1} successively constructed some saddle type nodal solutions  whose nodal domains are
of conical shapes demonstrating symmetric configurations by making use of
the dihedral groups and the polyhedral symmetric groups, respectively.  As a consequence, more general saddle solutions with prescribed Coxeter symmetries were constructed by a unified approach in  \cite{XiaZhang2020,Xia2021} for the critical and subcritical Choquard equations, respectively.

 Motivated by the above mentioned progresses on the nodal solutions of the Choquard equation \eqref{1.5}, we are led to a natural question: whether the saddle solutions for $s=1$ can still exist in the full range $s\in(0,1)$. In this aspect, some attempts have been made in \cite{Cui2020} by the first author of the current paper and the odd solution and the minimal nodal solution were obtained. It also turns out that the least energy on the Nehari nodal set is still the groundstate level for the case $p<2$. In the present paper, we will give an affirmative answer and construct saddle solutions for any prescribed Coxeter symmetric nodal structures.

 To state our results, we first introduce some notations.
 Let $G=\langle \mathcal{S} \rangle$ be a finite Coxeter group with $\mathcal{S}$ being its generating set and  $|\mathcal{S}|=k\in \{1, 2,\cdots, N\}$. Since $s^2=1$ for $s\in\mathcal{S}$, there exists a unique epimorphism $\phi:G\to \{\pm 1\}$ induced by $\phi(s)=-1$ for each $s\in \mathcal{S}$, see Lemma \ref{Lemma 2.2}.  A function $u:\bbR^N\to \bbR$ is said to be $G$-symmetry if it satisfies
$$g\circ u (x)=\phi(g)u(x), \quad \text{ for } g\in G,\quad x\in \bbR^N. $$
Here $\circ$ denotes  the group action that will be explained in Section \ref{section2}.
In what follows, we say $u\in H^s(\bbR^N)\setminus\{0\}$ is a $G$-saddle solution if $u$ solves the Choquard equation \eqref{eqfrc} with $G$-symmetry and $u$ is called $G$-groundstate if in addition $u$ minimizes the energy functional amongst all the $G$-saddle solutions.

Our main result can be stated as follows.
\begin{theorem}\label{Theorem 1}
Assume that $s\in(0,1)$, $N\geq3$, $\alpha\in (0,N)$ and $2\leq p<\frac{N+\alpha}{N-2s}$. For any Coxeter group $G$ with its rank $k$ satisfying $1\leq k\leq N$, the fractional Choquard equation \eqref{eqfrc} permits a $G$-saddle solution $u_{G}\in H^{s}(\mathbb{R}^{N})$ with exactly $\abs{G}$ nodal domains which turns to be a $G$-groundstate.
\end{theorem}

Although this seems to be a predictable result compared with the saddle solutions for the Choquard equation (see e.g., \cite{XW2,XW1,Xia2021}), it is worthwhile to remark that there are still some difficulties for the fractional Choquard equation and some new ideas are needed. On one hand, due to the appearance of the fractional Laplacian, solutions for the fractional Choquard equation \eqref{eqfrc} can not decay exponential anymore even for the case of $p\geq 2$, they turn out to be polynomial decay at infinity. This leads to a refine analysis in establishing the strict energy inequalities which shall play significant roles in restoring the compactness, see Proposition \ref{Proposition 3.1}; on the other hand, the signed property of the $G$-groundstate on the fundamental domain of $G$ can not be deduced directly by similar arguments as in \cite{Cui2020,MoS,XW1}. Actually, it will be proved by applying the strong maximum principle for the fractional Laplacian, see e.g., \cite[Corollary 4.12]{CS1}. To this end, we shall make use of the extension method for fractional Laplacian established by Caffarelli and Silvestre \cite{CS}, through which instead of the nonlocal operator we can study a local elliptic equation with a Neumann boundary condition in one dimension higher. As a result, our method takes advantage of finite Coxeter groups used in \cite{XiaZhang2020,Xia2021} and adapts it to the fractional Laplacian operator. Some ideas of \cite{Cui2020} are borrowed and our proofs turn more transparent.

This paper is organized as follows. Section \ref{section2} is devoted to introducing some preliminaries on the Coxeter groups and the variational framework space. The proofs of Theorem \ref{Theorem 1} will be completed in Section \ref{section3} and \ref{section4}. More precisely, the $G$-groundstate will be constructed inductively in Section \ref{section3} by employing the Lions' concentration-compactness principle. In the final section, we show the sign property of the $G$-groundstate that has exactly $\abs{G}$ nodal domains by using the extension method.

\section{Coxeter groups and variational framework}
\label{section2}
\subsection{Finite Coxeter groups}
 \quad\, In this section, we shall collect some basic results about Coxeter groups for the readers' convenience although it has been done in
 \cite[Sec. 2]{XiaZhang2020}. We refer to \cite{Da,Th} for further information. We first recall some notations and definitions. Let $H$ and $N$ are the subgroup of group $G$ and $N$ is normal. We write $G$ as $G=N\rtimes H$, if $H\cap N=\{1\}$ and $G=HN$. There are two basic concepts about group actions.

\begin{definition}\label{Definition 2.1}
 Suppose $G$ acts on $\mathbb{R}^{k}$ by homeomorphisms. The orbit of $x\in \mathbb{R}^{k}$ defined by $\mathcal{O}_{x}=\{gx|g\in G\}$.The isotropy subgroup of $x$ is $\cS_{x}=\{g\in G|gx=x\}$.
\end{definition}

For any set $M$,   $\abs{M}$ denotes the cardinality of $M$. The following Lagrange's theorem is well known.

\begin{lemma}\label{Lemma 2.1}
Assume that the group $G$ acts on $\mathbb{R}^{k}$. Then for any $x\in\mathbb{R}^{k}$, it holds
\begin{displaymath}\begin{aligned}
\abs{G}=\abs{\mathcal{O}_{x}}\abs{\cS_{x}}.
\end{aligned}\end{displaymath}
\end{lemma}

We now introduce the Coxeter symmetry. We begin with the definition of reflection. Let $V$ be a vector space. A linear reflection on $V$ is a linear automorphism $r: V\rightarrow V$ with $r^{2}=1$. And the wall $H_{r}$ is the set of midpoints of edges flipped by $r$. We then recall a formal definition of the Coxeter group.

\begin{definition}\label{Definition 2.2}
Let $I$ be a finite indexing set with cardinal number $|I|=k$ and let $\cS=\{s_{i}\}_{i\in I}$. Let $M=(M_{i,j})_{i,j\in I}$ be a $k\times k$ matrix such that $M_{i,j}=M_{j,i}\in\{1,2,\cdot\cdot\cdot,\infty\}$ for all $i,j\in I$ and $M_{i,j}=1$ if and only if  $i=j\in I$. Then $M$ is called a Coxeter matrix.  The associated Coxeter group $G_{\cS}$ is defined by the presentation
\begin{displaymath}\begin{aligned}
G=G_{\cS}=\langle\cS |(s_{i}s_{j})^{M_{i,j}}=1,\;\forall i,j\in I\rangle.
\end{aligned}\end{displaymath}
The pair $(G,\cS)$ is a Coxeter system and  $\cS$ is a Coxeter generating set of $G$. The cardinality $\abs{I}$ is usually called the rank of $(G,\cS)$.
\end{definition}
The following two results are fundamental for Coxeter groups.

\begin{lemma}\label{Lemma 2.2}
Let $(G,\cS)$  be a Coxeter system. Then there exists an epimorphism $\phi:G\rightarrow\{-1,1\}$ induced by $\phi(s)=-1$ for all $s\in \mathcal{S}$. Then for any $g\in G$, $\phi(g^{-1})=\phi^{-1}(g)=\phi(g)$.
\end{lemma}

\begin{lemma}\label{Lemma 2.3}
Let $(G,\cS)$  be a Coxeter system with rank $k$. Then there exists a faithful representation
 \begin{displaymath}\begin{aligned}
 \rho:G\rightarrow GL(\bbR^{k}).
 \end{aligned}\end{displaymath}
\end{lemma}

For any $\cT\subset \cS$, we can check $G_{\cT}$ is a Coxeter group. This leads to the irreducible Coxeter group.

\begin{definition}\label{Definition 2.3}
A Coxeter system $(G,\cS)$ is reducible if $\cS=\cS^{'}\cup \cS^{''}$ with $\cS^{'}$ and $\cS^{''}$ nonempty, such that
$(s_{i}s_{j})^{2}=1$ for all $s_{i}\in S^{'}$ and $s_{j}\in S^{''}$. It then follows that $G=G_{\cS^{'}}\times G_{\cS^{''}}$. A Coxeter system $(G,\cS)$ is irreducible if it is not reducible.
\end{definition}

For any finite irreducible Coxeter group, we have
\begin{lemma}\label{Lemma 2.4}
Let $(G,\cS)$ be a finite irreducible Coxeter system with rank $k$. Then for any $s\in \cS$, there exists a normal subgroup $N_{k_{s}}$ such that $G=N_{k_{s}}\rtimes \langle s\rangle $.
\end{lemma}

Finally, we recall the fundamental domain which will be used in later.
\begin{definition}\label{Definition 2.4}
Assume the Coxeter group $G$ acts on $\mathbb{R}^{k}$ by homeomorphisms. A fundamental domain is a closed, connected subset $\mathcal{D}$ of $\mathbb{R}^{k}$ such that $\mathcal{O}_{x}\cap\mathcal{D}\neq\emptyset$ for any $x\in \mathbb{R}^{k}$ and $\mathcal{O}_{x}\cap\mathcal{D}=\{x\}$ for any $x$ in the interior of $\mathcal{D}$.
\end{definition}
Note that $\cup_{g\in G}H_{g}$ separates $\mathbb{R}^{k}$ into $\abs{G}$ components. For any component, its closure is a fundamental domain of $G$. Write the $k-i$ dimensional facets of the fundamental domain $\mathcal{D}$ as $\partial^{i}\mathcal{D}$ for $i=0,1,\cdot\cdot\cdot,k-1$. We have

\begin{lemma}\label{Lemma 2.5}
Let $\mathcal{G}_{k}$ be the collection of the Coxeter group of rank $k$. Assume $\mathcal{D}$ be the fundamental domain of $G\in\mathcal{G}_{k}$. Then for any $x\in\partial^{i}\mathcal{D}$, $\mathcal{S}_{x}\in \mathcal{G}_{i}$ with $i=0,1,\cdot\cdot\cdot,k-1$.
\end{lemma}
\subsection{ Variational framework}
\quad\, We start with  the fractional Sobolev space. For a measurable function $u:\mathbb{R}^{N}\rightarrow\mathbb{R}$, the Gagliardo norm is defined by
\begin{displaymath}\begin{aligned}
\lceil u\rfloor_{s}=\Big(\int_{\mathbb{R}^{N}}\int_{\mathbb R^N}\frac{\abs{u(x)-u(y)}^{2}}{\abs{x-y}^{N+2s}}\dif x\dif y\Big)^{\frac{1}{2}}.
\end{aligned}\end{displaymath}
And by \cite[Propositions 3.4 and 3.6]{DPV}, we have
\begin{equation}\label{2.1}\begin{aligned}
2\int_{\mathbb{R}^{N}}\abs{\xi}^{2s}\abs{\mathcal{F}(u)(\xi)}^{2}\dif \xi=2\|(-\Delta)^{s/2}u\|^{2}_{L^{2}(\mathbb{R}^{N})} =C_{N,s}\lceil u\rfloor^{2}_{s}.
\end{aligned}\end{equation}
  Then the fractional Sobolev space can be defined by
\begin{displaymath}
\begin{aligned}
H^{s}(\mathbb{R}^{N})=\{u\in L^{2}(\mathbb{R}^{N})|\lceil u\rfloor_{s}<+\infty\},
\end{aligned}\end{displaymath}
 endowed with the norm
\begin{displaymath}\begin{aligned}
\|u\|_{2,s}^{2}=\|u\|_{L^{2}(\mathbb{R}^{N})}^{2}+\lceil u\rfloor^{2}_{s}.
\end{aligned}\end{displaymath}
Up to a constant, we may assume $C_{N,s}=2$, i.e., we define for $u\in H^s(\mathbb R^N)$ that
$$
\|(-\Delta)^{s/2}u\|_{L^{2}(\mathbb{R}^{N})}=\lceil u\rfloor_s.
$$
The following fractional embedding  theorem has been proved, see \cite[Theorem 6.5]{DPV} for instance.
\begin{lemma}\label{Lemma 2.6} Let $s\in(0,1)$ and $N>2s$, then there exists
a positive constant $C=C(N,s)>0$ such that, for any $u\in H^{s}(\mathbb{R}^{N})$
 \begin{displaymath}\begin{aligned}
\|u\|_{L^{2^{*}_{s}}(\mathbb{R}^{N})}\leq C\|u\|_{2,s},
\end{aligned}\end{displaymath}
with $2^{*}_{s}=\frac{2N}{N-2s}$. Moreover the embedding $H^{s}(\mathbb{R}^{N})\hookrightarrow L^{r}(\mathbb{R}^{N})$ is continuous for any $r\in[2,2^{*}_{s}]$.
\end{lemma}

We shall study the fractional Choquard equation \eqref{eqfrc} with the help of  Coxeter's symmetries. Let $G\in \mathcal{G}_k$ be a Coxeter group with its rank $1\leq k\leq N$. For any $g\in G$,
the group action on  $u\in H^{s}(\mathbb{R}^{N})$ is defined by
\begin{displaymath}\begin{aligned}
g\circ u=u(g^{-1}x), \quad \text{ with }\quad gx=\diag(g,1_{N-k})x.
\end{aligned}\end{displaymath}
We shall seek  saddle solutions for the fractional Choquard equation \eqref{eqfrc} in the following subspace
\begin{displaymath}
\begin{aligned}
H_{G}^{s}(\mathbb{R}^{N})=\{u\in H^{s}(\mathbb{R}^{N})| \, g\circ u(x)=\phi(g)u(x),\, \forall g\in G\},
\end{aligned}\end{displaymath}
where $\phi$ is the unique epimorphism induced by $\phi(s)=-1$ for $s\in \mathcal{S}$, see Lemma \ref{Lemma 2.2}.

The fractional Choquard equation is variational in nature, and it is easy to see weak solutions of \eqref{eqfrc} correspond to  critical points of the following action functional
\begin{displaymath}\begin{aligned}
I(u)=\frac{1}{2}\int_{\mathbb{R}^{N}}\abs{(-\Delta)^{s/2}u}^{2}+\abs{u}^{2}\dif x-
\frac{1}{2p}\int_{\mathbb{R}^{N}}\int_{\mathbb R^N}\frac{A_{\alpha}\abs{u(y)}^{p}\abs{u(x)}^{p}}{\abs{x-y}^{N-\alpha}}\dif x\dif y.
\end{aligned}\end{displaymath}
From the Hardy--Littlewood--Sobolev inequality
and the fractional embedding theorem,  we see that the energy functional $I:H_G^s(\bbR^N)\to \bbR$ is well defined and belongs to $C^{1}$ if $p$ satisfies \eqref{1.3}.
It is known that the positive groundstate can be describe as
\begin{equation*}\label{1.4}
\begin{aligned}
c_{0}=\inf_{v\in\mathcal{N}}I(v),
\end{aligned}\end{equation*}
where $\mathcal{N}$ is the Nehari constraint
\begin{displaymath}\begin{aligned}
\mathcal{N}=\{u\in H^{s}(\mathbb{R}^{N})\setminus\{0\} : \langle I'(u),u\rangle=0\}.
\end{aligned}\end{displaymath}

We shall consider the minimization problem
\begin{displaymath}\begin{aligned}
c_{G}=\inf\limits_{u\in\mathcal{N}_{G}}I(u),
\end{aligned}\end{displaymath}
with the $G$-Nehari manifold being defined by
\begin{displaymath}\begin{aligned}
\mathcal{N}_{G}=\mathcal{N}\cap H_{G}^{s}(\mathbb{R}^{N}).
\end{aligned}\end{displaymath}
On the other hand, the above minimal energy has a mountain pass type description. In fact, we can conclude as \cite[Theorem 4.2]{W} that
\begin{equation}\label{2.2}\begin{aligned}
c_G=\inf\limits_{\gamma\in\Gamma}\max_{t\in[0,1]}I(\gamma(t))
\end{aligned}\end{equation}
where the paths set is
\begin{displaymath}\begin{aligned}
\Gamma=\{\gamma\in C([0,1];H_{G}^{s}(\mathbb{R}^{N}))|\gamma(0)=0,\;I(\gamma(1))<0\}.
\end{aligned}\end{displaymath}

The qualitative properties including regularity and decay behaviours of solutions for the fractional Choquard equation \eqref{eqfrc} have been obtained in \cite[Theorems 3.2-3.3]{DSS} and it turns to be of importance in seeking saddle solutions. We recall and adopt it as follows.
\begin{lemma}\label{Lemma 3.1} Assume that $N\geq 3$, $\alpha\in (0,N)$ and $2\leq p<\frac{N+\alpha}{N-2s}$. If $v\in H^{s}(\mathbb{R}^{N})$ with $s\in(0,1)$ is a solution of \eqref{eqfrc}, then $v\in L^{1}(\mathbb{R}^{N})\cap C^{\beta}(\bbR^N)$ for some  $\beta\in(0,1)$. 
Moreover there exists $C>0$ such that  for all $x\in\mathbb{R}^{N}$,
$$\abs{v(x)}\leq C(1+\abs{x}^{2})^{-(N+2s)/2}.$$
\end{lemma}
\section{Existence of saddle solutions}
\label{section3}

\quad\,  We begin with an energy estimate.
\begin{proposition}\label{Proposition 3.1}
Let $s\in(0,1)$, $N\geq 3$, $\alpha\in(0,N)$ and $2\leq p<\frac{N+\alpha}{N-2s}$. For any $G\in \mathcal{G}_{k}$ with its rank satisfying $1\leq k\leq N$, there exists $c^*_{G}$ such that
\begin{align*}
0<c_{G}<c^*_{G}&=\min\{\abs{\mathcal{O}_{x}}c_{\mathcal{S}_{x}}\mid x\in \partial^{k-1}\mathcal{D}\}&\\
&<\min\{\abs{\mathcal{O}_{x}}c_{\mathcal{S}_{x}}\mid x\in \partial^{k-i}\mathcal{D},i=2,\cdots,k\}.
\end{align*}
\end{proposition}

\noindent{\bf{Proof}}. Since $\mathcal{N}_{G}\subset\mathcal{N}$ then $0<c_{0}\leq c_{G}$. The remaining inequalities will be completed by the method of induction. Assume $k=1$. By the classification results of the finite Coxeter group (see Theorem 6.9.1 in \cite{Da}), $G=\rm{A}_{1}$, the cyclic group with order 2, and its fundamental domain is, up to a rotation, $\mathbb{R}^{N}_{+}$. In this case, the inequality to be proved reads as $c_{G}=c_{\rm{odd}}<2c_{0}$, which is exactly the result in \cite{Cui2020}.

We now suppose that Proposition \ref{Proposition 3.1} holds true for $k-1$. Then $H$-saddle solutions for the fractional Choquard equation \eqref{eqfrc} do exist for any $H\leq G$ with $H\in \cup_{i<k}\mathcal{G}_i$, this can be proved by repeating our subsequent procedures. In particular, each $H$-saddle solution we obtained is an $H$-groundstate.
Fix $q\in\partial^{k-1}\mathcal{D}$ such that $\abs{q}=\sqrt{q\cdot q}=1$. Without loss of generality, we assume that $q$ is
perpendicular to $H_{s_1}$.
By Lemma \ref{Lemma 2.5}, we see that $\mathcal{S}_{q}\in\mathcal{G}_{k-1}$. The inductive assumption implies that the fractional Choquard equation \eqref{eqfrc} has a saddle solution    $u_{\cS_q}$   with $ \mathcal{S}_{q}$-symmetry such that
\begin{equation}\label{3.1}\begin{aligned}
c_{ \mathcal{S}_{q}}=I(u_{\mathcal{S}_{q}})=\inf_{ \cN_{ \mathcal{S}_q}} I ,\; \text{ and } \; g\circ u_{\mathcal{S}_{q}}=\phi(g)u_{\mathcal{S}_{q}},\;\forall g\in  \mathcal{S}_{q},
\end{aligned}\end{equation}
where $\phi:  \cS_q\rightarrow\{-1,1\}$ is the restriction of the unique epimorphism induced by $\phi(s)=-1$ for $s\in \mathcal{S}$.
Moreover, from Lemma \ref{Lemma 3.1}, it follows that
\begin{equation}\label{3.2}\begin{aligned}
u_{\mathcal{S}_{q}}\in L^{1}(\mathbb{R}^{N})\cap L^{\infty}(\mathbb{R}^{N}),\text{ and }\limsup_{\abs{x}\to+\infty}\abs{x}^{N+2s}\abs{u_{\mathcal{S}_{q}}}<+\infty.
\end{aligned}\end{equation}

As a result of Definition \ref{Definition 2.3}, a finite Coxeter group always can be reduced to the direct sum of some finite irreducible Coxeter groups, see \cite[Theorem 6.9.1]{Da}. Without loss of generality, we  write $G=G_{k_{1}}\times G_{r}$ where $G_{k_{1}}\in \mathcal{G}_{k_1}$ is a finite irreducible Coxeter group such that $s_{1}\in G_{k_{1}}$ and $G_r\in \mathcal{G}_{k-k_1}$. If $G$ is irreducible,
then $G_r=\{1\}$.
For the case $k_{1}=1$, we  assume that $H_{s_{1}}=\partial\mathbb{R}^{N}_{+}$. We define a function $u_{R}:\mathbb{R}^{N}\rightarrow\mathbb{R}$ for each $x=(x', x^{N})\in\mathbb{R}^{N-1}\times\mathbb{R}$ by
\begin{displaymath}\begin{aligned}
u_{R}=(\xi_{R}u_{\mathcal{S}_{q}})(x', x^{N}-3R)-(\xi_{R}u_{\mathcal{S}_{q}})(x', -x^{N}-3R),
\end{aligned}\end{displaymath}
where the function $\xi_{R}(x)=\xi(x/R)$ and $\xi(x)\in C^{\infty}(\mathbb{R}^{N})$ is radial and satisfies $\xi(x)=1$ if $\abs{x}\leq1$, $\xi(x)=0$ if $\abs{x}\geq2$ and $0\leq\xi\leq1$ on $\mathbb{R}^{N}$. In this situation, we can easily check $G_r=\langle \mathcal{S}_{q}\rangle$. Since $s_{1}g=gs_1$ for any $g\in G_{r}$. We can treat $G_{r}$ as a group acting on $\mathbb R^{N-1}$ so that $g\circ u_{R}=\phi(g)u_{R}$ for any $g\in G_{r}$.  Hence by $s_{1}\circ u_{R}=-u_{R}$, we then conclude that $u_{R}\in H^{s}_{G}(\mathbb{R}^{N})$.

For the case $k_{1}\geq2$, by Lemma \ref{Lemma 2.4}, there exists $N_{k_1}\lhd G$ such that $G_{k_{1}}=N_{k_{1}}\rtimes \langle s_{1}\rangle$. Let $N_{k}=N_{k_{1}}\times G_{r}$. We now defined a function $u_{R}$ by
\begin{displaymath}\begin{aligned}
u_{R}=\frac{1}{\abs{ \mathcal{S}_{q}} }\sum_{g\in N_{k}}(\xi_{R}u_{\mathcal{S}_{q}})(g^{-1}s_{1}x-l_{G}Rq)-(\xi_{R}u_{\mathcal{S}_{q}})(g^{-1}x-l_{G}Rq),
\end{aligned}\end{displaymath}
where $l_{G}$ is a constant such that $l_{G}\min\limits_{x\neq y\in\mathcal{O}_{q}}\abs{x-y}\geq6$. In this case, we can check that $ \mathcal{S}_{q}=G_{S_{k_{1r}}}\times G_{r}$, where $S_{k_{1r}}=\cS\setminus\{s_1\}=\{s_{2},\cdot\cdot\cdot,s_{k_{1}}\}$. We now show that $u_{R}\in H_{G}^{s}(\mathbb{R}^{N})$.
From the definition of $u_R$, we first observe that $s_{1}\circ u_{R}=-u_{R}=\phi(s_{1})u_{R}$. For $s\in S_{k_{1r}}$, since $s\in G_{k_{1}}=N_{k_{1}}\rtimes \langle s_{1}\rangle$, there exists unique $n_{s},h_{s}$ such that $s=n_{s}s_{1}=s_{1}h_{s}$. Hence, we can deduce that
\begin{displaymath}\begin{aligned}
s\circ u_{R}&=\frac{1}{\abs{ {\mathcal{S}_{q}}}}\sum_{g\in N_{k}}(\xi_{R}u_{\mathcal{S}_{q}})(g^{-1}s_{1}s^{-1}x-l_{G}Rq)-(\xi_{R}u_{\mathcal{S}_{q}})(g^{-1}s^{-1}x-l_{G}Rq)\\
&=\frac{1}{\abs{ {\mathcal{S}_{q}}}}\sum_{g\in N_{k}}(\xi_{R}u_{\mathcal{S}_{q}})(g^{-1}n^{-1}_{s}x-l_{G}Rq)
-(\xi_{R}u_{\mathcal{S}_{q}})(g^{-1}h^{-1}_{s}s_{1}x-l_{G}Rq)\\
&=-\frac{1}{\abs{ {\mathcal{S}_{q}}}}\sum_{g\in N_{k}}(\xi_{R}u_{\mathcal{S}_{q}})(g^{-1}s_{1}x-l_{G}Rq)-(\xi_{R}u_{\mathcal{S}_{q}})(g^{-1}x-l_{G}Rq)=-u_{R}.\\
\end{aligned}\end{displaymath}
Note that $gh=hg$ for $g\in G_{r}$ and $h\in G_{k_{1}}$. We deduce that
\begin{displaymath}\begin{aligned}
g_{r}\circ u_{R}&=\frac{1}{\abs{ {\mathcal{S}_{q}}}}\sum_{g\in N_{k}}(\xi_{R}u_{\mathcal{S}_{q}})(g^{-1}s_{1}g_{r}^{-1}x-l_{G}Rq)
-(\xi_{R}u_{\mathcal{S}_{q}})(g^{-1}g_{r}^{-1}x-l_{G}Rq)\\
&=\phi(g_{r})u_{R}.
\end{aligned}\end{displaymath}
The conclusion $u_{R}\in H_{G}^{s}(\mathbb{R}^{N})$ follows immediately by combining $s_1\circ u_R=-u_R$.

Take $t_{R}>0$ such that $\langle I'(t_{R}u_{R}),t_{R}u_{R}\rangle=0$. We then deduce that,
\begin{equation}\label{energyesti}
\begin{aligned}
I(t_{R}u_{R})=\Big(\frac{1}{2}-\frac{1}{2p}\Big)\frac{\displaystyle \Big(\int_{\mathbb{R}^{N}}\abs{(-\Delta)^{s/2}u_{R}}^{2}
+\abs{u_{R}}^{2}\dif x\Big)^{\frac{p}{p-1}}}
{\displaystyle \Big(\int_{\mathbb{R}^{N}}(K_{\alpha}\ast\abs{u_{R}}^{p})\abs{u_{R}}^{p}\dif x\Big)^{\frac{1}{p-1}}}.
\end{aligned}\end{equation}
Therefore, $c_{G}<c^*_{G}$ follows once we establish that for some $R>0$,

\begin{align*}\label{3.3}
\frac{\Big(\displaystyle \int_{\mathbb{R}^{N}}\abs{u_{R}}^{2}+\abs{(-\Delta)^{s/2}u_{R}}^{2}\dif x\Big)^{\frac{p}{p-1}}}
{\Big(\displaystyle \int_{\mathbb{R}^{N}}(K_{\alpha}\ast\abs{u_{R}}^{p})\abs{u_{R}}^{p}\dif x\Big)^{\frac{1}{p-1}}}
<
\abs{\mathcal{O}_{q}}\frac{\displaystyle \Big(\int_{\mathbb{R}^{N}}\abs{u_{\mathcal{S}_{q}}}^{2}+\abs{(-\Delta)^{s/2}u_{\mathcal{S}_{q}}}^{2}\dif x\Big)^{\frac{p}{p-1}}}
{\Big(\displaystyle \int_{\mathbb{R}^{N}}(K_{\alpha}\ast\abs{u_{\mathcal{S}_{q}}}^{p})\abs{u_{\mathcal{S}_{q}}}^{p}\dif x\Big)^{\frac{1}{p-1}}}.
\end{align*}

Let $u_{R}^{g}=(\xi_{R} u_{\mathcal{S}_{q}})(g^{-1}x-l_{G}Rq)$ and $D_{g}=\overline{\supp  u_{R}^{g}} $. Here the notation $\supp v$ denotes the support of the function $v$. Direct calculus give us that
\begin{equation}\label{3.4}\begin{aligned}
\|(-\Delta)^{s/2}u_{R}\|_{L^{2}(\mathbb{R}^{N})}^{2}=&\frac{1}{\abs{ {\mathcal{S}_{q}}}}\sum\limits_{g\in G}
\|(-\Delta)^{s/2}u_{R}^{g}\|^{2}_{L^{2}(\mathbb{R}^{N})}\\
&+\frac{1}{\abs{ {\mathcal{S}_{q}}}}\sum\limits_{g\neq h\in G}
\int_{\mathbb{R}^{N}}(-\Delta)^{s/2}u_{R}^{g}(-\Delta)^{s/2}u_{R}^{h}\dif x.
\end{aligned}\end{equation}
Thanks to \eqref{2.1}, we deduce that
\begin{equation}\label{3.5}\begin{aligned}
\|(-\Delta)^{s/2}(\xi_{R}u_{\mathcal{S}_{q}})\|_{L^{2}(\mathbb{R}^{N})}^{2}
=\|(-\Delta)^{s/2}u_{R}^{g}\|_{L^{2}(\mathbb{R}^{N})}^{2},\;\forall g\in G.
\end{aligned}\end{equation}
By the definition of $\xi$ and the choice of $l_{G}$, we have
\begin{equation}\label{3.6}
\begin{aligned}
2R\leq \min\limits_{x\in D_{g},y\in D_{h}}\abs{x-y}\leq(2l_{G}+4)R,\; \forall g\neq h.
\end{aligned}\end{equation}
Then, by combining \eqref{2.1}, \eqref{3.4}-\eqref{3.6} and the fact $u_{\mathcal{S}_{q}}\in L^{1}(\mathbb{R}^{N})$, we deduce that
\begin{equation}\label{3.7}\begin{aligned}
\|(-\Delta)^{s/2}u_{R}\|_{L^{2}(\mathbb{R}^{N})}^{2}
\leq\abs{\mathcal{O}_{q}}\|(-\Delta)^{s/2}(\xi_{R}u_{\mathcal{S}_{q}})\|_{L^{2}(\mathbb{R}^{N})}^{2}+
\frac{C}{R^{N+2s}}\|u_{\mathcal{S}_{q}}\|^{2}_{L^{1}(\mathbb{R}^{N})}.
\end{aligned}\end{equation}

Let
\begin{displaymath}\begin{aligned}
u_{R}^{0}(x,y)=(\xi_{R}u_{\mathcal{S}_{q}})(x)-(\xi_{R}u_{\mathcal{S}_{q}})(y),\;x,y\in\mathbb{R}^{N}.
\end{aligned}\end{displaymath}
We see from \eqref{2.1} that
\begin{displaymath}
\begin{aligned}
2C_{N,s}^{-1}\|(-\Delta)^{s/2}(\xi_{R}u_{\mathcal{S}_{q}})\|_{L^{2}(\mathbb{R}^{N})}^{2}
=\sum^{3}_{i=1}\sum^{3}_{j=1}\int_{E_{i}}\int_{E_{j}}\frac{\abs{u_{R}^{0}(x,y)}^{2}}{\abs{x-y}^{N+2s}}\dif x\dif y=
\sum^{3}_{i=1}\sum^{3}_{j=1}F_{ij},
\end{aligned}\end{displaymath}
where $E_{1}=B_{R}$, $E_{2}=B_{2R}\setminus B_{R}$, $E_{3}=\mathbb{R}^{N}\setminus B_{2R}$. By the definition of $\xi_{R}$,
\begin{displaymath}\begin{aligned}
F_{11}=\int_{B_{R}}\int_{B_{R}}
\frac{\abs{u_{\mathcal{S}_{q}}(x)-u_{\mathcal{S}_{q}}(y)}^{2}}{\abs{x-y}^{N+2s}}\dif x\dif y.
\end{aligned}\end{displaymath}
By the decay behaviour of $u_{\mathcal{S}_{q}}$, we deduce that
\begin{equation} \label{decayproperty}
\begin{aligned}
\int_{\mathbb{R}^{N}\setminus B_{R}}\abs{u_{\mathcal{S}_{q}}(y)}\dif y\leq C\int_{\mathbb{R}^{N}\setminus B_{R}}\frac{1}{\abs{y}^{N+2s}}\dif y= {C}{R^{-2s}}.
\end{aligned}\end{equation}
Take $r,t>1$ such that
$1/r+1/t+(N+2s-2)/N=2$, so that $N/r+N/t=N+2-2s$. Note that $\abs{\xi'(t)}\leq 2$.
We deduce by the Hardy--Littlewood--Sobolev inequality and \eqref{decayproperty} that
\begin{equation}\label{importantestimate1}
\begin{aligned}
\int_{B_{2R}\setminus B_{R}}\int_{ B_{R}}&\frac{\abs{\xi_{R}(x)-\xi_R(y)}^2\abs{u_{\mathcal{S}_{q}}(y)}^2}{\abs{x-y}^{N+2s}}\dif x\dif y\\
&
\leq \frac{C}{R^2}\int_{B_{2R}\setminus B_{R}}\int_{ B_{R}}\frac{\abs{x-y}^2\abs{u_{\mathcal{S}_{q}}(y)}^2}{\abs{x-y}^{N+2s}}\dif x\dif y\\
&
\leq \frac{C}{R^2}\Big(\int_{B_R} 1\dif x\Big)^{1/r}\Big(\int_{B_{2R}\setminus B_R}\abs{u_{\mathcal{S}_q}}^{2t}\dif y\Big)^{1/t}\\
&
\leq \frac{C}{R^2} \frac{R^{N/r}}{R^{2(N+2s)-N/t}}\leq \frac{C}{R^{N+6s}}=o(R^{-8s}).
\end{aligned}
\end{equation}
Similarly, we have
\begin{align}\label{importantestimate2}
&\int_{B_{2R}\setminus B_{R}}\int_{ B_{2R}\setminus B_R}\frac{\abs{\xi_{R}(x)-\xi_R(y)}^2\abs{u_{\mathcal{S}_{q}}(y)}^2}{\abs{x-y}^{N+2s}}\dif x\dif y
=o(R^{-8s}),
\end{align}
and
\begin{align}\label{importantestimate3}
&\int_{\bbR^N\setminus B_{2R}}\int_{ B_{2R}\setminus B_R}\frac{\abs{\xi_{R}(x)-\xi_{R}(y)}^2\abs{u_{\mathcal{S}_{q}}(y)}}{\abs{x-y}^{N+2s}}\dif x\dif y
=O(R^{-4s}).
\end{align}
We observe that
\begin{displaymath}
\begin{aligned}
F_{21}&=\int_{B_{2R}\setminus B_{R}}\int_{B_{R}}\frac{\abs{u_{\mathcal{S}_{q}}(x)-(\xi_{R}u_{\mathcal{S}_{q}})(y)}^{2}}
{\abs{x-y}^{N+2s}}\dif x\dif y\\
&=\int_{B_{2R}\setminus B_{R}}\int_{B_{R}}\frac{\abs{u_{\mathcal{S}_{q}}(x)-u_{\mathcal{S}_{q}}(y)}^{2}}
{\abs{x-y}^{N+2s}}\dif x\dif y\\
&\quad +\int_{B_{2R}\setminus B_{R}}\int_{ B_{R}}\frac{\abs{\xi_{R}(x)-\xi_{R }(y)|^{2}|u_{\mathcal{S}_{q}}(y)}^{2}}
{\abs{x-y}^{N+2s}}\dif x\dif y\\
&\quad+2\int_{B_{2R}\setminus B_{R}}\int_{B_{R}}\frac{(u_{\mathcal{S}_{q}}(x)-u_{\mathcal{S}_{q}}(y))(\xi_{R}(x)-\xi_{R}(y))
u_{\mathcal{S}_{q}}(y)}{\abs{x-y}^{N+2s}}\dif x\dif y.
\end{aligned}
\end{displaymath}
Note that  $u_{\mathcal{S}_q}\in H^s(\bbR^N)$ has  bounded norm.
By combining \eqref{importantestimate1} and the Cauchy--Schwartz inequality, we deduce that
\begin{align*}
&\int_{B_{2R}\setminus B_{R}}\int_{B_{R}}\frac{(u_{\mathcal{S}_{q}}(x)-u_{\mathcal{S}_{q}}(y))(\xi_{R}(x)-\xi_{R}(y))
u_{\mathcal{S}_{q}}(y)}{\abs{x-y}^{N+2s}}\dif x\dif y\\
&\leq  \Big(\int_{B_{2R}\setminus B_{R}}\int_{B_{R}} \frac{\abs{u_{\mathcal{S}_{q}}(x)-u_{\mathcal{S}_{q}}(y)}^2}{\abs{x-y}^{N+2s}}\dif x\dif y\Big)^{1/2}\\
&\quad \times  \Big(\int_{B_{2R}\setminus B_{R}}\int_{B_{R}} \frac{\abs{\xi_{R}(x)-\xi_{R}(y)}^2\abs{
u_{\mathcal{S}_{q}}(y)}^2}{\abs{x-y}^{N+2s}}\dif x\dif y\Big)^{1/2}\\
&\leq C\Big(\int_{B_{2R}\setminus B_{R}}\int_{B_{R}} \frac{\abs{\xi_{R}(x)-\xi_{R}(y)}^2\abs{
u_{\mathcal{S}_{q}}(y)}^2}{\abs{x-y}^{N+2s}}\dif x\dif y\Big)^{1/2}=o(R^{-4s}).
\end{align*}
It then follows that
\begin{displaymath}
\begin{aligned}
F_{21} 
&\leq\int_{B_{2R}\setminus B_{R}}\int_{B_{R}}\frac{\abs{u_{\mathcal{S}_{q}}(x)-u_{\mathcal{S}_{q}}(y)}^{2}}
{\abs{x-y}^{N+2s}}\dif x\dif y+o(R^{-4s}).
\end{aligned}\end{displaymath}
In a similar manner, we can conclude  by \eqref{importantestimate2}  that
\begin{displaymath}
\begin{aligned}
F_{22}&=\int_{B_{2R}\setminus B_{R}}\int_{B_{2R}\setminus B_{R}}\frac{\abs{(\xi_{R}u_{\mathcal{S}_{q}})(x)-(\xi_{R}u_{\mathcal{S}_{q}})(y)}^{2}}
{\abs{x-y}^{N+2s}}\dif x\dif y\\
&=\int_{B_{2R}\setminus B_{R}}\int_{B_{2R}\setminus B_{R}}\frac{\xi_{R}(x)\xi_{R}(y)\abs{u_{\mathcal{S}_{q}}(x)-u_{\mathcal{S}_{q}}(y)}^{2}}
{\abs{x-y}^{N+2s}}\dif x\dif y\\
&\quad+ \int_{B_{2R}\setminus B_{R}}\int_{B_{2R}\setminus B_{R}}\frac{\abs{\xi_R(x)-\xi_R(y)}^2\abs{u_{\mathcal{S}_{q}}(y)}^{2}}{\abs{x-y}^{N+2s}}\dif x\dif y\\
&\quad+ \int_{B_{2R}\setminus B_{R}}\int_{B_{2R}\setminus B_{R}}\frac{\xi_R(x)(\xi_R(x)-\xi_R(y))(u^2_{\mathcal{S}_{q}}(x)-u^2_{\mathcal{S}_{q}}(y))}
{\abs{x-y}^{N+2s}}\dif x\dif y\\
&\leq\int_{B_{2R}\setminus B_{R}}\int_{B_{2R}\setminus B_{R}}\frac{\abs{u_{\mathcal{S}_{q}}(x)-u_{\mathcal{S}_{q}}(y)}^{2}}
{\abs{x-y}^{N+2s}}\dif x\dif y+o(R^{-4s}).
\end{aligned}\end{displaymath}
Recall that $\xi_{R}\equiv0$ in $\mathbb{R}^{N}\setminus B_{2R}$. By \eqref{importantestimate3} and $u_{\mathcal{S}_{q}}\in L^{\infty}(\mathbb{R}^{N})$ , we have
\begin{displaymath}
\begin{aligned}
F_{32}&=\int_{\mathbb{R}^{N}\setminus B_{2R}}\int_{B_{2R}\setminus B_{R}}\frac{\abs{(\xi_{R}u_{\mathcal{S}_{q}})(x)}^{2}}
{\abs{x-y}^{N+2s}}\dif x\dif y\\
&=\int_{\mathbb{R}^{N}\setminus B_{2R}}\int_{B_{2R}\setminus B_{R}}
\frac{\xi_{R}^{2}(x)\abs{u_{\mathcal{S}_{q}}(x)-u_{\mathcal{S}_{q}}(y)}^{2}}{\abs{x-y}^{N+2s}}\dif x\dif y\\
&\quad+\int_{ \bbR^N\setminus B_{2R}}\int_{B_{2R}\setminus B_{R}}
\frac{\abs{\xi_{R}(x)-\xi_{R}(y)}^{2}(2u_{\mathcal{S}_{q}}(x)-u_{\mathcal{S}_{q}}(y))u_{\mathcal{S}_{q}}(y)}
{\abs{x-y}^{N+2s}}\dif x\dif y\\
&\leq \int_{\mathbb{R}^{N}\setminus B_{2R}}\int_{B_{2R}\setminus B_{R}}
\frac{\abs{u_{\mathcal{S}_{q}}(x)-u_{\mathcal{S}_{q}}(y)}^{2}}{\abs{x-y}^{N+2s}}\dif x\dif y+O(R^{-4s}).
\end{aligned}\end{displaymath}
Note that $\abs{x-y}\geq R$ for $x\in B_R$ and $y\in \bbR^N\setminus B_{2R}$. We conclude that
\begin{displaymath}
\begin{aligned}
F_{31}&=\int_{\mathbb{R}^{N}\setminus B_{2R}}\int_{B_{R}}\frac{\abs{u_{\mathcal{S}_{q}}(x)}^{2}}
{\abs{x-y}^{N+2s}}\dif x\dif y\\
&=\int_{\mathbb{R}^{N}\setminus B_{2R}}\int_{ B_{R}}\frac{\abs{u_{\mathcal{S}_{q}}(x)-u_{\mathcal{S}_{q}}(y)}^{2}}
{\abs{x-y}^{N+2s}}\dif x\dif y\\
&\quad +\int_{\mathbb{R}^{N}\setminus B_{2R}}\int_{B_{R}}\frac{(2u_{\mathcal{S}_{q}}(x)-u_{\mathcal{S}_{q}}(y))u_{\mathcal{S}_{q}}(y)}
{\abs{x-y}^{N+2s}}\dif x\dif y \\
&\leq\int_{\mathbb{R}^{N}\setminus B_{2R}}\int_{ B_{R}}\frac{\abs{u_{\mathcal{S}_{q}}(x)-u_{\mathcal{S}_{q}}(y)}^{2}}
{\abs{x-y}^{N+2s}}\dif x\dif y+\frac{CR^{N}}{R^{N+2s}}\int_{\mathbb{R}^{N}\setminus B_{2R}}\abs{u_{\mathcal{S}_{q}}(y)}\dif y\\
&=\int_{\mathbb{R}^{N}\setminus B_{2R}}\int_{ B_{R}}\frac{\abs{u_{\mathcal{S}_{q}}(x)-u_{\mathcal{S}_{q}}(y)}^{2}}
{\abs{x-y}^{N+2s}}\dif x\dif y+O(R^{-4s}).\\
\end{aligned}\end{displaymath}
Observing that
\begin{displaymath}
\begin{aligned}
&F_{12}=F_{21}\leq \int_{ B_{R}}\int_{B_{2R}\setminus B_{R}}\frac{|u_{\mathcal{S}_{R}}(x)-u_{\mathcal{S}_{R}}(y)|^{2}}{\abs{x-y}^{N+2s}}\dif x\dif y+o(R^{-4s}),\\
&F_{13}=F_{31}\leq \int_{B_{R}}\int_{\mathbb{R}^{N}\setminus B_{2R}}\frac{\abs{u_{\mathcal{S}_{q}}(x)-u_{\mathcal{S}_{q}}(y)}^{2}}{\abs{x-y}^{N+2s}}\dif x\dif y+O(R^{-4s}),\\
&F_{23}=F_{32}\leq \int_{B_{2R}\setminus B_{R}}\int_{\mathbb{R}^{N}\setminus B_{2R}}\frac{\abs{u_{\mathcal{S}_{q}}(x)-u_{\mathcal{S}_{q}}(y)}^{2}}{\abs{x-y}^{N+2s}}\dif x\dif y+O(R^{-4s}),
\end{aligned}\end{displaymath}
we finally obtain by noting $F_{33}=0$ that
\begin{equation*}
\begin{aligned}
\|(-\Delta)^{s/2}(\xi_{R}u_{\mathcal{S}_{q}})\|_{L^{2}(\mathbb{R}^{N})}^{2}=\frac{C_{N,s}}{2}\sum_{1\leq i,j\leq 3}F_{ij}\leq \lceil u_{\mathcal{S}_{q}}\rfloor^{2}_{s}+O(R^{-4s}).
\end{aligned}
\end{equation*}
This, together with \eqref{2.1} and \eqref{3.7}, implies that
\begin{equation*}\label{3.12}
\begin{aligned}
\|(-\Delta)^{s/2}u_{R}\|^{2}_{L^{2}(\mathbb{R}^{N})}\leq \abs{\mathcal{O}_{q}}\|(-\Delta)^{s/2}u_{\mathcal{S}_{q}}\|^{2}_{L^{2}(\mathbb{R}^{N})}+O(R^{-4s}).
\end{aligned}\end{equation*}
 By the construction of $u_{R}$, we easily obtain
\begin{equation*}
\begin{aligned}
\int_{\mathbb{R}^{N}}\abs{u_{R}}^{2}\dif x=\abs{\mathcal{O}_{q}}\int_{B_{2R}}\abs{\xi_{R}u_{\mathcal{S}_{q}}}^{2}\dif x\leq \abs{\mathcal{O}_{q}}\int_{\mathbb{R}^{N}}\abs{u_{\mathcal{S}_{q}}}^{2}\dif x.
\end{aligned}\end{equation*}
Note that $(N+\alpha)/(N-2s)>2$. we deduce that
\begin{equation}\label{3.14}
\begin{aligned}
\|u_R\|_{2,s}^{2}\leq \abs{\mathcal{O}_{q}}\int_{\mathbb{R}^{N}}\abs{u_{\mathcal{S}_{q}}}^{2}+\abs{(-\Delta)^{s/2}u_{\mathcal{S}_{q}}}^{2}\dif x
+o(R^{\alpha-N}).
\end{aligned}\end{equation}

We now deal with the denominator of \eqref{energyesti}. By \eqref{3.6}, we have
\begin{displaymath}
\begin{aligned}
&\int_{\mathbb{R}^{N}}(K_{\alpha}\ast\abs{u_{R}}^{p})\abs{u_{R}}^{p}\dif x\\
&=\abs{\mathcal{O}_{q}}\int_{\mathbb{R}^{N}}(K_{\alpha}\ast\abs{\xi_{R}u_{\mathcal{S}_{q}}}^{p})
\abs{\xi_{R}u_{\mathcal{S}_{q}}}^{p}\dif x
+\frac{1}{\abs{ {\mathcal{S}_{q}}}}\sum_{g\neq h}\int_{\mathbb{R}^{N}}(K_{\alpha}\ast\abs{u_{R}^{g}}^{p})\abs{u_{R}^{h}}^{p}\dif x\\
&\geq\abs{\mathcal{O}_{q}}\int_{\mathbb{R}^{N}}(K_{\alpha}\ast\abs{\xi_{R}u_{\mathcal{S}_{q}}}^{p})
\abs{\xi_{R}u_{\mathcal{S}_{q}}}^{p}\dif x
+\frac{C}{R^{N-\alpha}}\big(\int_{B_{R}}\abs{u_{\mathcal{S}_{q}}}^{p}\dif x\big)^{2}.
\end{aligned}\end{displaymath}
For the first term, we have
\begin{displaymath}
\begin{aligned}
&\int_{\mathbb{R}^{N}}(K_{\alpha}\ast\abs{\xi_{R}u_{\mathcal{S}_{q}}}^{p})\abs{\xi_{R}u_{\mathcal{S}_{q}}}^{p}\dif x\\
=&\int_{\mathbb{R}^{N}}(K_{\alpha}\ast\abs{u_{\mathcal{S}_{q}}}^{p})\abs{u_{\mathcal{S}_{q}}}^{p}\dif x
-2\int_{\mathbb{R}^{N}}(K_{\alpha}\ast\abs{u_{\mathcal{S}_{q}}}^{p})(1-\xi^{p}_{R}) \abs{u_{\mathcal{S}_{q}}}^{p}\dif x\\
&\quad+\int_{\mathbb{R}^{N}}(K_{\alpha}\ast(1-\xi^{p}_{R})\abs{u_{\mathcal{S}_{q}}}^{p})
(1-\xi^{p}_{R})\abs{u_{\mathcal{S}_{q}}}^{p}\dif x\\
\geq&\int_{\mathbb{R}^{N}}(K_{\alpha}\ast\abs{u_{\mathcal{S}_{q}}}^{p})\abs{u_{\mathcal{S}_{q}}}^{p}\dif x
-2\int_{\mathbb{R}^{N}}(K_{\alpha}\ast\abs{u_{\mathcal{S}_{q}}}^{p})(1-\xi^{p}_{R})\abs{u_{\mathcal{S}_{q}}}^{p}\dif x.
\end{aligned}\end{displaymath}
By the decay properties of $u_{\mathcal{S}_{q}}$, we see that $u_{\mathcal{S}_{q}}\in L^r(\bbR^N)$ for $r>\frac{N}{N+2s}$. Moreover,
for every $r>\frac{\alpha}{N+2s}$, we conclude similarly as in \cite[Proposition  4.1]{XiaZhang2020} that $K_\alpha\ast|u_{\mathcal S_q}|^r \in L^\infty(\bbR^N)$. It then follows that
there exists positive constant $C>0$ such that
 \begin{displaymath}
\begin{aligned}
\limsup_{\abs{x}\rightarrow +\infty}\frac{K_{\alpha}\ast\abs{u_{\mathcal{S}_{q}}}^{p}}{K_{\alpha}(x)}\leq C,
\end{aligned}
\end{displaymath}
which leads that
\begin{displaymath}
\begin{aligned}
2\int_{\mathbb{R}^{N}}(K_{\alpha}\ast\abs{u_{\mathcal{S}_{q}}}^{p})(1-\xi^{p}_{R})\abs{u_{\mathcal{S}_{q}}}^{p}\dif x
\leq C\int_{\mathbb{R}^{N}\setminus B_{R}}\frac{\abs{u_{\mathcal{S}_{q}}}^{p}}{\abs{x}^{N-\alpha}}\dif x.
\end{aligned}\end{displaymath}
Therefore, by the asymptotic properties of $u_{\mathcal{S}_{q}}$ again, we deduce that
\begin{align}\label{3.15}
&\int_{\mathbb{R}^{N}}(K_{\alpha}\ast\abs{u_{R}}^{p})\abs{u_{R}}^{p}\dif x &\nonumber\\
&\geq\abs{\mathcal{O}_{q}}\int_{\mathbb{R}^{N}}(K_{\alpha}\ast\abs{u_{\mathcal{S}_{q}}}^{p})\abs{u_{\mathcal{S}_{q}}}^{p}\dif x
+\frac{C}{R^{N-\alpha}}\Big(\int_{B_{R}}\abs{u_{\mathcal{S}_{q}}}^{p}\dif x\Big)^{2}-C\int_{\mathbb{R}^{N}\setminus B_{R}}\frac{\abs{u_{\mathcal{S}_{q}}}^{p}}{\abs{x}^{N-\alpha}}\dif x&\nonumber\\
&\geq\abs{\mathcal{O}_{q}}\int_{\mathbb{R}^{N}}(K_{\alpha}\ast\abs{u_{\mathcal{S}_{q}}}^{p})\abs{u_{\mathcal{S}_{q}}}^{p}\dif x
+\frac{C}{R^{N-\alpha}}\Big(\int_{B_{R}}\abs{u_{\mathcal{S}_{q}}}^{p}\dif x\Big)^{2}+o( R^{\alpha-N}).
\end{align}
We then conclude by combining \eqref{3.14}, \eqref{3.15} and the definition of $c_G$ that
$$c_G \leq I(t_Ru_R)=\abs{\mathcal{O}_q}c_{ {\cS_q}}\Big(1-CR^{\alpha-N}+o\big(R^{\alpha-N}\big)\Big)<\abs{\mathcal{O}_q}c_{ {\cS_q}}.$$

We now in position to give the existence of  saddle solution with $G$-symmetry.

\noindent{\bf Proof of Theorem \ref{Theorem 1}} Thanks to \eqref{2.2}, by the general minimax principle (see e.g., \cite[Theorem 2.8]{W}), we can obtain a sequence $\{u_{n}\}_{n\geq1}\subset H^{s}_{G}(\mathbb{R}^{N})$ such that
\begin{equation}\label{3.16}
\begin{aligned}
I(u_{n})\rightarrow c_{G},\quad I'(u_{n})\rightarrow0\quad in \quad(H^{s}_{G}(\mathbb{R}^{N}))^{\ast} \quad as\quad n\rightarrow\infty.
\end{aligned}\end{equation}
It is easy to verify that $\{u_{n}\}_{n\geq1}\subset H^{s}_{G}(\mathbb{R}^{N})$ is bounded, since
\begin{displaymath}
\begin{aligned}
c_{G}+o(1)\|u_{n}\|_{2,s}
=I(u_{n})-\frac{1}{2p}\langle I'(u_{n}),u_{n}\rangle=(\frac{1}{2}-\frac{1}{2p})\|u_n\|_{2,s}^2.
\\
\end{aligned}\end{displaymath}
We shall claim that there exist $T>0$ and $\{a_{n}\}_{n\geq1}\subset\{0\}\times\mathbb{R}^{N-k}$ such that
\begin{equation}\label{3.17}\begin{aligned}
\liminf\limits_{n\rightarrow\infty}\int_{B_{T}(a_{n})}\abs{u_{n}}^{\frac{2Np}{N+\alpha}}\dif x>0.
\end{aligned}\end{equation}
Up to translations and a subsequence, we may assume $u_{n}$ converges weakly to some function $u_{G}\in H_{G}^{s}(\mathbb{R}^{N})\setminus\{0\}$. A standard procedure (e.g., \cite{GS}, see also \cite{Cui2020,XW1}) implies that $I'(u_{G})=0$ in $(H_G^s(\bbR^N))^*$ and $I(u_{G})=c_{G}$. We then apply the symmetric criticality principle \cite[Theorem 1.28]{W} to conclude that $u_G$ is a critical point of $I$ in $H^s(\bbR^N)$.

To prove \eqref{3.17}, we first assert that
\begin{equation}\label{3.18}
\begin{aligned}
\lim\limits_{n\rightarrow\infty}\int_{\mathbb{R}^{N}}
\abs{u_{n}(x)}^{\frac{2Np}{N+\alpha}}\dif x=\Lambda\in(0,+\infty).
\end{aligned}\end{equation}
Indeed, by observing  $\{u_{n}\}_{n\geq1}$ in $H^s(\mathbb R^N)$, we have $\Lambda<+\infty$ according to  the fractional embedding theorem (Lemma \ref{Lemma 2.6}). And if $\Lambda=0$, we would deduce a contradiction by the Hardy--Littlewood--Sobolev inequality that
\begin{align*}0<c_{G}=(\frac{1}{2}-\frac{1}{2p})\lim\limits_{n\rightarrow\infty}\int_{\mathbb{R}^{N}}\int_{\mathbb R^N}
\frac{A_{\alpha}\abs{u_{n}(x)}^{p}\abs{u_{n}(y)}^{p}}{\abs{x-y}^{N-\alpha}}\dif x\dif y=0.\end{align*}

We first observe that  the sequence $\{u_n\}_{n\geq 1}$ is non-vanishing, e.g., for any $r>0$, we have
$\lim_{n\to\infty}\sup_{y\in\mathbb R^N}\int_{B_r(y)}|u_n|^{\frac{2Np}{N+\alpha}}\dif x>0$.
Otherwise, we conclude by \cite[Lemma 2.3]{DSS}  that $\{u_n\}_{n\geq 1}$  converges to zero strongly in $L^r(\mathbb R^N)$ for any $r\in (2,\frac{2N}{N-2s})$, which contradicts \eqref{3.18} since $\Lambda>0$.
According to Lions' concentration compactness lemma \cite{Li2}, we shall consider the remaining two cases.

{\bf Compactness}: there exists $\{x_{n}\}_{n\geq1}\subset \mathbb{R}^{N}$, such that for any $\varepsilon>0$, there exists $R_{0}>0$ such that
\begin{align*}
\liminf\limits_{n\rightarrow\infty}\int_{B_{R_{0}}(x_{n})}\abs{u_{n}(x)}^{\frac{2Np}{N+\alpha}}\dif x
\geq\Lambda-\varepsilon.
\end{align*}
In this case, we can easily check that for $\{x_n\}_{n\geq 1}$, there exists some $M_1>0$ such that
\begin{equation}\label{3.19}
\begin{aligned}
\abs{\Prj_{k}(x_{n})}\leq M_{1}.
\end{aligned}\end{equation}
Here $\Prj_{k}:\mathbb R^N\to\mathbb R^k, (x^1,\cdots,x^N)\mapsto (x^{1},\cdots,x^{k})$ denotes the projection. Otherwise, for every $g\in G\setminus\{1\}$ and for largely $n$,  we have that
$B_{R_0}(x_{n})\cap gB_{R_0}(x_{n})=\emptyset$. We then deduce a contradiction by the symmetric setting that
\begin{align*}
\liminf\limits_{n\rightarrow\infty}\int_{\mathbb{R}^{N}}\abs{u_{n}(x)}^{\frac{2Np}{N+\alpha}}\dif x
\geq&\liminf\limits_{n\rightarrow\infty}\int_{B_{R_0}(x_{n})}\abs{u_{n}(x)}^{\frac{2Np}{N+\alpha}}\dif x\\
&+\liminf\limits_{n\rightarrow\infty}\int_{gB_{R_0}(x_{n})}\abs{u_{n}(x)}^{\frac{2Np}{N+\alpha}}\dif x
\geq2\Lambda-2\varepsilon>\Lambda.
\end{align*}

 {\bf  Dichotomy}: there exist $\beta\in(0,1)$ and  $\{y_{n}\}_{n\geq1}\subset\mathbb{R}^{N}$, such that for any $\varepsilon>0$, there exists $R_{2}>0$ such that for all $r_{1}\geq R_{2}$ and $r_{2}\geq R_{2}$,
\begin{displaymath}
\begin{aligned}
\limsup\limits_{n\rightarrow\infty}\Big|\int_{B_{r_{1}}(y_{n})}\abs{u_{n}(x)}^{\frac{2Np}{N+\alpha}}\dif x-\beta\Lambda\Big|
+\Big|\int_{\mathbb{R}^{N}\setminus B_{r_{2}}(y_{n})}\abs{u_{n}(x)}^{\frac{2Np}{N+\alpha}}\dif x-(1-\beta)\Lambda\Big|\leq\varepsilon.
\end{aligned}
\end{displaymath}
In this situation, we will also show that
\begin{equation}\label{3.20}
\begin{aligned}
\abs{\Prj_{k}(y_{n})}\leq M_{2},\;\forall n,
\end{aligned}\end{equation}
for some appropriate constant $M_{2}>0$. Without loss of generality, we may assume $\{y_{n}\}_{n\geq 1}\subset \mathcal{D}$, the fundamental domain of the Coxeter group $G$.
If there exists $R_1>0$ such that
$
\{y_n\}_{n\geq 1}\subset B_{R_1}\cap \mathcal{D}
$,
 the conclusion \eqref{3.20} can be deduced immediately  since $\beta\Lambda>0$. For the subsequent proofs, we assume by contradiction that $\abs{\Prj_{k}(y_{n})}\rightarrow +\infty$, then by a diagonal process and up to a subsequence, we can choose  $\varepsilon_{n}\rightarrow0$,  and $r_{n}'=4r_{n}\rightarrow +\infty$,
 such that
$$\Big|\int_{B_{r_n}(y_n)}\abs{u_n}^{\frac{2Np}{N+\alpha}}-\beta\Lambda\Big|
+\Big|\int_{B^c_{r'_n}(y_n)}\abs{u_n}^{\frac{2Np}{N+\alpha}}-(1-\beta)\Lambda\Big|\leq \varepsilon_n.
$$
As a result, several cases occur due to the variant of the asymptotic distance between the sequence $\{y_n\}_{n\geq 1}$ and the boundary $\partial\mathcal D$.\\
\textbf{Case 1}: up to a subsequence, there exists some $j\in\{1,\cdots,k-1\}$ such that $$\lim_{n\to\infty}\dist(y_n, \partial^j \mathcal D)=0 \text{ and }
\lim_{n\to\infty}\dist(y_n, \partial^\ell \mathcal D)=0 \text{ for all } \ell>j.$$
In this situation, we fix a unit vector $q\in \partial^j\mathcal D\cap \partial B_1$ and, up to a subsequence, we may assume that $y_n=\nu_nq$ with $\nu_n\to+\infty$ for each subcases such that
$\abs{x-y}\geq 4r_n$ for $x\in B_{3r_n}(z_1)$ and $y\in B_{3r_n}(z_2)$ where $z_1\neq z_2\in Gy_n$.
For any $g\in G$, we denote $B_{r_{n}}^{g}=B_{r_{n}}(gy_{n})$ and $\mathcal{B}_{r_{n}}=\bigcup\limits_{g\in G}B_{r_{n}}^{g}$. Then we have
\begin{equation}\label{3.21}\begin{aligned}
\int_{\mathcal{B}_{4r_{n}}\setminus \mathcal{B}_{r_{n}}}\abs{u_{n}}^{\frac{2Np}{N+\alpha}}\dif x\leq\abs{\mathcal{O}_{y_{n}}}\varepsilon_{n}.
\end{aligned}\end{equation}
We now prove that
\begin{equation}\label{3.22}
\begin{aligned}
\int_{\mathcal{B}_{3r_{n}}\setminus \mathcal{B}_{2r_{n}}}
\int_{\mathbb{R}^{N}}\frac{\abs{u_{n}(x)-u_{n}(y)}^{2}}{\abs{x-y}^{N+2s}}\dif y\dif x+\int_{\mathcal{B}_{3r_{n}}\setminus \mathcal{B}_{2r_{n}}}\abs{u_{n}}^{2}\dif x=o(1).
\end{aligned}\end{equation}
Take $\psi_{r}(x)=\psi(x/r)$ where $\psi\in C^{\infty}(\mathbb{R}^{N})$ satisfies that $\psi(x)=0$ for $\abs{x}\leq1$ or $\abs{x}\geq4$, $\psi(x)=1$ for $2\leq\abs{x}\leq3$ and $0\leq\psi(x)=\psi(\abs{x})\leq1$. For any $g\in G$, let $\psi_{n}^{g}(x)=\psi_{r_{n}}(x-gy_{n})$ and
\begin{displaymath}\begin{aligned}
\Psi_{n}(x)=\sum\limits_{z\in\mathcal{O}_{y_{n}}}\psi_{r_{n}}(x-z)
=\frac{1}{\abs{\mathcal{S}_{y_{n}}}}\sum\limits_{g\in G}\psi_{n}^{g}(x).
\end{aligned}\end{displaymath}
Observe that for any $h\in G$,
\begin{displaymath}\begin{aligned}
h\circ(\Psi_{n}u_{n})(x)&=\frac{1}{\abs{ {\mathcal{S}_{y_{n}}}}}\sum\limits_{g\in G}\psi_{r_{n}}(h^{-1}x-h^{-1}hgy_{n})u_{n}(h^{-1}x)\\
&=\frac{1}{\abs{ {\mathcal{S}_{y_{n}}}}}\sum\limits_{g\in G}\psi_{r_{n}}(x-h gy_{n})\phi(h)u_{n}(x)=\phi(h)(\Psi_{n}u_{n})(x),
\end{aligned}\end{displaymath}
we infer that $\Psi_{n}u_{n}\in H^{s}_{G}(\mathbb{R}^{N})$.  Since $I'(u_{n})\rightarrow0$ in $ (H^{s}_{G}(\mathbb{R}^{N}))^{\ast}$, we have $\langle I'(u_{n}),\Psi_{n}u_{n}\rangle=o(1)$. And a direct computation shows that
 \begin{align*}
\int_{\mathbb{R}^{N}}\int_{\mathbb R^N}&\frac{\Psi_{n}(x)\abs{u_{n}(x)-u_{n}(y)}^{2}}{\abs{x-y}^{N+2s}}\dif x\dif y
+\int_{\mathbb{R}^{N}}\Psi_{n} \abs{u_{n}}^{2}\dif x\\
=&\langle I'(u_{n}),\Psi_{n}u_{n}\rangle+\int_{\mathcal{B}_{4r_{n}}\setminus \mathcal{B}_{r_{n}}}
(K_{\alpha}\ast\abs{u_{n}}^{p})\abs{u_{n}(x)}^{p}\Psi_{n}(x)\dif x\\
&\quad-\int_{\mathbb{R}^{N}}\int_{\mathbb R^N}\frac{(\Psi_{n}(x)-\Psi_{n}(y))(u_{n}(x)-u_{n}(y))u_{n}(y)}{\abs{x-y}^{N+2s}}\dif x\dif y.
\end{align*}
Then by the Cauchy--Schwartz inequality, the Hardy--Littlewood--Sobolev inequality and \eqref{3.21}, we obtain that
\begin{equation}\label{3.23}\begin{aligned}
\int_{\mathbb{R}^{N}}\int_{\mathbb R^N}&\frac{\Psi_{n}(x)\abs{u_{n}(x)-u_{n}(y)}^{2}}{\abs{x-y}^{N+2s}}\dif x\dif y
+\int_{\mathbb{R}^{N}}\Psi_{n} \abs{u_{n}}^{2}\dif x\\
\leq&\lceil u_{n}\rfloor_{s}
\Big(\int_{\mathbb{R}^{N}}\int_{\mathbb R^N}\frac{\abs{\Psi_{n}(x)-\Psi_{n}(y)}^{2}\abs{u_{n}(y)}^{2}}{\abs{x-y}^{N+2s}}\dif x\dif y\Big)^{1/2}+o(1).
\end{aligned}\end{equation}
By the choice of $\Psi_{n}$, we readily check that
\begin{align}\label{infinitesimalestimate1}
\int_{\mathbb{R}^{N}}&\int_{\mathbb R^N}\frac{\abs{\Psi_{n}(x)-\Psi_{n}(y)}^{2}\abs{u_{n}(y)}^{2}}{\abs{x-y}^{N+2s}}\dif x\dif y \nonumber\\
&=\frac{1}{\abs{ {\mathcal{S}_{y_{n}}}}}\sum\limits_{g\in G}\int_{\mathbb{R}^{N}}\int_{\mathbb R^N}
\frac{\abs{\psi_{n}^{g}(x)-\psi^{g}_{n}(y)}^{2}\abs{u_{n}(y)}^{2}}{\abs{x-y}^{N+2s}}\dif x\dif y & \\
&\quad -\frac{1}{\abs{ {\mathcal{S}_{y_{n}}}}}\sum\limits_{g\neq h\in G}\int_{\mathbb{R}^{N}}\int_{\mathbb R^N}
\frac{\psi_{n}^{g}(x)\psi^{h}_{n}(y)\abs{u_{n}(y)}^{2}}{\abs{x-y}^{N+2s}}\dif x\dif y\nonumber\\
&\quad -\frac{1}{\abs{ {\mathcal{S}_{y_{n}}}}}\sum\limits_{g\neq h\in G}\int_{\mathbb{R}^{N}}\int_{\mathbb R^N}
\frac{\psi_{n}^{g}(y)\psi^{h}_{n}(x)\abs{u_{n}(y)}^{2}}{\abs{x-y}^{N+2s}}\dif x\dif y\nonumber\\
&=\frac{1}{\abs{ {\mathcal{S}_{y_{n}}}}}\sum\limits_{g\in G}\int_{\mathbb{R}^{N}}\int_{\mathbb R^N}
\frac{\abs{\psi_{n}^{g}(x)-\psi^{g}_{n}(y)}^{2}\abs{u_{n}(y)}^{2}}{\abs{x-y}^{N+2s}}\dif x\dif y
+C\|u_{n}\|^{2}_{L^{2}(\mathbb{R}^{N})}r_{n}^{-2s}.\nonumber
\end{align}
On one hand, we have
\begin{align}\label{infinitesimalestimate2}
 &\int_{\mathbb{R}^{N}}\int_{\mathbb R^N}
\frac{\abs{\psi_{n}^{g}(x)-\psi^{g}_{n}(y)}^{2}\abs{u_{n}(y)}^{2}}{\abs{x-y}^{N+2s}}\dif x\dif y \\
&\leq \int_{\mathbb{R}^{N}}\int_{B_{4r_{n}}(gy_{n})}
\frac{\abs{\psi_{n}^{g}(x)-\psi^{g}_{n}(y)}^{2}\abs{u_{n}(y)}^{2}}{\abs{x-y}^{N+2s}}\dif x\dif y\nonumber\\
&\quad +  \int_{B_{4r_{n}}(gy_{n})}\int_{B_{5r_{n}}(gy_{n})}
\frac{\abs{\psi_{n}^{g}(x)-\psi^{g}_{n}(y)}^{2}\abs{u_{n}(y)}^{2}}{\abs{x-y}^{N+2s}}\dif x\dif y\nonumber\\
&\quad +\int_{B_{4r_{n}}(gy_{n})}\int_{\mathbb{R}^{N}\setminus B_{5r_{n}}(gy_{n})}
\frac{\abs{\psi_{n}^{g}(x)-\psi^{g}_{n}(y)}^{2}\abs{u_{n}(y)}^{2}}{\abs{x-y}^{N+2s}}\dif x\dif y.\nonumber
\end{align}
By using the Hardy--Littlewood--Sobolev inequality, we deduce for any  $l\in(0,+\infty)$ that,
\begin{align}\label{infinitesimalestimate3}
 \int_{\mathbb{R}^{N}}\int_{B_{lr_{n}}(gy_{n})}
&\frac{\abs{\psi_{n}^{g}(x)-\psi^{g}_{n}(y)}^{2}\abs{u_{n}(y)}^{2}}{\abs{x-y}^{N+2s}}\dif x\dif y\nonumber\\
&\leq \frac{C}{r_n^2} \int_{\mathbb{R}^{N}}\int_{B_{lr_{n}}(gy_{n})}
\frac{\abs{x-y}^{2}\abs{u_{n}(y)}^{2}}{\abs{x-y}^{N+2s}}\dif x\dif y\nonumber\\
&\leq \Big(\int_{B_{lr_n}(gy_n)} 1 \dif x\Big)^{\frac{1}{r}}\Big(\int_{\bbR^N}\abs{u_n(y)}^{2t}\dif y\Big)^{\frac{1}{t}}\leq \frac{C}{r_n^{2-N/r}}=o(1).
\end{align}
Here we take $r,t>1$ such that $1/r+1/{t}+(N+2s-2)/{N}=2$ and $t<N/(N-2s)$  so that $2-{N}/{r}=2s-N+N/{t}>0$.
On the other hand,
\begin{align}\label{infinitesimalestimate4}
\int_{B_{4r_{n}}(gy_{n})}&\int_{\mathbb{R}^{N}\setminus B_{5r_{n}}(gy_{n})}
\frac{\abs{\psi_{n}^{g}(x)-\psi^{g}_{n}(y)}^{2}\abs{u_{n}(y)}^{2}}{\abs{x-y}^{N+2s}}\dif x\dif y\nonumber\\
&\leq C\int_{B_{4r_n}(gy_n)} \abs{u_n(y)}^2 \dif y\int_{\bbR^N\setminus B_{r_n}(y)}\frac{1}{\abs{x-y}^{N+2s}}\dif x\nonumber\\
&\leq C\int_{r_n}^{+\infty}\frac{\rho^{N-1}}{\rho^{N+2s}}\dif \rho=\frac{C}{r_n^{2s}}.
\end{align}
Inserting these estimates \eqref{infinitesimalestimate2}--\eqref{infinitesimalestimate4} into \eqref{infinitesimalestimate1}, we deduce
\begin{equation}\label{3.24}\begin{aligned}
\lim_{n\to\infty}\int_{\mathbb{R}^{N}}\int_{\mathbb R^N}\frac{\abs{\Psi_{n}(x)-\Psi_{n}(y)}^{2}\abs{u_{n}(y)}^{2}}{\abs{x-y}^{N+2s}}\dif x\dif y=0.
\end{aligned}\end{equation}
This, together with \eqref{3.23} yields the conclusion \eqref{3.22}.

We now take another radial cut-off function $\eta\in C^{\infty}(\mathbb{R}^{N})$ such that $\eta(x)=1$ if $x\in B_{2}$, $\eta(x)=0$ if $|x|\geq3$ and $0\leq\eta(x)=\eta(|x|)\leq1$. Let $\eta_{r}(x)=\eta(x/r)$. For any
$1\leq i\leq|\mathcal{O}_{y_{n}}|$ and $z^{i}\in\mathcal{O}_{y_{n}}$, we denote $v^{i}_{n}=\eta_{r_{n}}(x-z^{i})u_{n}$. Let $v_{n}=\sum_iv_{n}^{i}$
and $w_{n}=u_{n}-v_{n}$. It is clear that $v_{n}, w_{n}\in H^{s}_{G}(\mathbb{R}^{N})$.  Particularly, we have
$v_{n}^{i}\in H_{G_{\mathcal{S}_{q}}}^{s}(\mathbb{R}^{N})$ with $G_{\mathcal {S}_q}\in \mathcal{G}_j$ for every $i$. Indeed, if $z^{i}=\nu_{n}g_iq$,
 then for every $s\in \mathcal{S}_q$, we deduce
$$g_isg_i^{-1}\diamond v_n^i =\eta(\abs{g_isg_i^{-1}x-z_i}/r_n)u_n(g_isg_i^{-1}x)=\psi(g_isg_i^{-1})v_n^i=\psi(s) v_n^i,$$
so that $v_n^i\in H_{g_iG_{\mathcal{S}_{q}}g_i^{-1}}^{s}(\mathbb{R}^{N})$. The conclusion follows since
the group $g_iG_{\mathcal{S}_q}g_i^{-1}$ is isomorphic to $G_{\mathcal{S}_q}$  for every $i$.

We now claim that
\begin{equation}\label{3.26}
\begin{aligned}
I(u_{n})=\sum_{i=1}^{\abs{\mathcal{O}_{y_{n}}}}I(v^{i}_{n})+I(w_{n})+o(1),
\end{aligned}\end{equation}
and
\begin{equation}\label{3.27}\begin{aligned}
\langle I'(v^{i}_{n}),v^{i}_{n}\rangle=o(1),\; \langle I'(w_{n}),w_{n}\rangle=o(1).
\end{aligned}\end{equation}
Assume \eqref{3.26} and \eqref{3.27} hold. By  Lemma \ref{Lemma 2.6}  and the Hardy-Littlewood-Sobolev inequality, we have
\begin{displaymath}\begin{aligned}
\lim\limits_{n\rightarrow\infty}\int_{\mathbb{R}^{N}}\abs{(-\Delta)^{s/2}v^{i}_{n}}^{2}+\abs{v^{i}_{n}}^{2}\dif x
=\lim\limits_{n\rightarrow\infty}\int_{\mathbb{R}^{N}}\int_{\mathbb R^N}
\frac{A_{\alpha}\abs{v^{i}_{n}(x)}^{p}\abs{v^{i}_{n}(y)}^{p}}{\abs{x-y}^{N-\alpha}}\dif x\dif y>0.
\end{aligned}\end{displaymath}
Then for each $1\leq i\leq \mathcal{O}_{y_n}$, there exists a unique $t^{i}_{n}\in(0,+\infty)$ such that $t^{i}_{n}v^{i}_{n}\in \mathcal{N}_{\mathcal{S}_{q}}$ satisfying that $\lim_{n\to\infty} t_n^i=1$.
Hence by \eqref{3.26} and \eqref{3.27}, we  deduce that
\begin{displaymath}\begin{aligned}
c_{G}&=\lim\limits_{n\rightarrow\infty}I(u_{n})
=\lim\limits_{n\rightarrow\infty}\sum_{i=1}^{\abs{\mathcal{O}_{y_{n}}}}I(t^{i}_{n}v^{i}_{n})+\lim\limits_{n\rightarrow\infty}I(w_{n})\\
&\geq\lim\limits_{n\rightarrow\infty}\sum_{i=1}^{\abs{\mathcal{O}_{y_{n}}}}
I(t^{i}_{n}v^{i}_{n})+\Big(\frac{1}{2}-\frac{1}{2p}\Big)
\lim\limits_{n\rightarrow\infty}\int_{\mathbb{R}^{N}}\abs{(-\Delta)^{s/2}w_{n}}^{2}
+\abs{w_{n}}^{2}\dif x\\
&\geq \abs{\mathcal{O}_{y_{n}}}c_{\mathcal{S}_{y_{n}}},
\end{aligned}\end{displaymath}
which is in contradiction to Proposition \ref{Proposition 3.1} that $c_{G}<\abs{\mathcal{O}_q} c_{\mathcal{S}_q}$.
This contradiction implies \eqref{3.20}.

 We turn now to proofs of the claims \eqref{3.26} and \eqref{3.27}. Indeed, we readily verify that
\begin{equation}\label{3.28}
\begin{aligned}
I(u_{n})=&\sum_{i=1}^{\abs{\mathcal{O}_{y_{n}}}}I(v^{i}_{n})+I(w_{n})+\int_{\mathbb{R}^{N}}(-\Delta)^{s/2}v_{n}(-\Delta)^{s/2}w_{n}\dif x
+\int_{\mathbb{R}^{N}}v_{n}w_{n}\dif x\\
&+\sum_{i\neq m}\int_{\mathbb{R}^{N}}(-\Delta)^{s/2}v^{i}_{n}(-\Delta)^{s/2}v^{m}_{n}\dif x
-\frac{1}{2p}\int_{\mathbb{R}^{N}}\int_{\mathbb{R}^{N}}F_n(x,y)\dif x\dif y\\
&-\sum_{i\neq m}\frac{1}{2p}\int_{\mathbb{R}^{N}}\int_{\mathbb{R}^{N}}
\frac{A_{\alpha}\abs{v^{i}_{n}(x)}^{p}\abs{v^{m}_{n}(y)}^{p}}{\abs{x-y}^{N-\alpha}}\dif x\dif y,
\end{aligned}\end{equation}
where $F_n:\bbR^N\times\bbR^N\to \bbR$ is
\begin{displaymath}\begin{aligned}
F_n(x,y)=A_{\alpha}\Big(\frac{\abs{u_{n}(x)}^{p}\abs{u_{n}(y)}^{p}}{\abs{x-y}^{N-\alpha}}-\frac{\abs{v_{n}(x)}^{p}\abs{v_{n}(y)}^{p}}{\abs{x-y}^{N-\alpha}}
-\frac{\abs{w_{n}(x)}^{p}\abs{w_{n}(y)}^{p}}{\abs{x-y}^{N-\alpha}}\Big).
\end{aligned}\end{displaymath}
Note that $\abs{x-y}\geq4r_{n}$ for $x\in \supp v^{i}_{n}$ and $y\in \supp v^{m}_{n}$ with $i\neq m$.
By the Cauchy--Schwartz inequality, we obtain for $i\neq m$ that

\begin{equation}\label{3.29}
\begin{aligned}
\int_{\mathbb{R}^{N}}& (-\Delta)^{s/2}v^{i}_{n}(-\Delta)^{s/2}v^{m}_{n}\dif x\\
&\leq  \int_{\mathbb{R}^{N}}\int_{\bbR^N}\frac{\abs{v^{i}_{n}(x)-v^{i}_{n}(y)}\abs{v^{m}_{n}(x)-v^{m}_{n}(y)}}{\abs{x-y}^{N+2s}}\dif x\dif y \\
&\leq 2\int_{B^{i}_{3r_{n}}}\int_{B^{m}_{3r_{n}}}
\frac{\abs{v^{i}_{n}(y)v^{m}_{n}(x)}}{\abs{x-y}^{N+2s}}\dif x\dif y\\
&\leq \frac{C r_{n}^{N}}{r_{n}^{N+2s}}\Big(\int_{B^{i}_{3r_{n}}}\abs{u_{n}(y)}^{2}\dif y\Big)^{1/2} \Big(\int_{B^{m}_{3r_{n}}}
\abs{u_{n}(x)}^{2}\dif x\Big)^{1/2}\\
&\leq\frac{C}{r_{n}^{2s}}\int_{\mathbb{R}^{N}}\abs{u_{n}}^{2}\dif x=o(1).
\end{aligned}\end{equation}
By the fractional embedding theorem, we deduce for $i\neq m$ that
\begin{equation}\label{3.30}
\begin{aligned}
&\int_{\mathbb{R}^{N}}\int_{\mathbb{R}^{N}}
\frac{\abs{v^{i}_{n}(x)}^{p}\abs{v^{m}_{n}(y)}^{p}}{\abs{x-y}^{N-\alpha}}\dif x\dif y
\leq\frac{C\|u_n\|^{2p}_{L^{p}(\mathbb{R}^{N})}}{r_{n}^{N-\alpha}}
\leq\frac{C\|u_{n}\|^{2p}_{2,s}}{r_{n}^{N-\alpha}}=o(1).
\end{aligned}\end{equation}
Let $\eta_{n}=\sum_{i=1}^{\abs{\mathcal{O}_{y_{n}}}}\eta_{r_{n}}(x-z^{i})$. By the choice of $\eta$, we have
\begin{displaymath}\begin{aligned}
& \int_{\mathbb{R}^{N}}(-\Delta)^{s/2}v_{n}(-\Delta)^{s/2}w_{n}\dif x +\int_{\mathbb{R}^{N}} v_{n}w_{n} \dif x\\
=&\int_{\mathbb{R}^{N}}\int_{\bbR^N}\frac{\eta_{n}(x)(1-\eta_{n}(x))\abs{u_{n}(x)-u_{n}(y)}^{2}}
{\abs{x-y}^{N+2s}} + \frac{\abs{\eta_{n}(x)-\eta_{n}(y)}^{2}u^{2}_{n}(y)}
{\abs{x-y}^{N+2s}}\dif x\dif y\\
&+  \int_{\mathbb{R}^{N}}\int_{\bbR^N}\frac{(1-\eta_{n}(x))(u_{n}(x)-u_{n}(y))(\eta_{n}(x)-\eta_{n}(y))u_{n}(y)}
{\abs{x-y}^{N+2s}}\dif x\dif y\\
&+ \int_{\mathbb{R}^{N}}\int_{\bbR^N}\frac{\eta_{n}(x)(u_{n}(x)-u_{n}(y))(\eta_{n}(x)-\eta_{n}(y))u_{n}(y)}
{\abs{x-y}^{N+2s}}\dif x\dif y +\int_{\mathbb{R}^{N}}v_{n}w_{n}\dif x\\
\leq&\int_{\mathcal{B}_{3r_{n}}\setminus \mathcal{B}_{2r_{n}}}
\int_{\mathbb{R}^{N}}\frac{\abs{u_{n}(x)-u_{n}(y)}^{2}}{\abs{x-y}^{N+2s}}\dif y\dif x+\int_{\mathbb{R}^{N}}\int_{\bbR^N}\frac{\abs{\eta_{n}(x)-\eta_{n}(y)}^{2}\abs{u_{n}(y)}^2}
{\abs{x-y}^{N+2s}}\dif x\dif y\\
&+\lceil u_{n}\rfloor_{s}\Big(\int_{\mathbb{R}^{N}}\int_{\bbR^N}\frac{|\eta_{n}(x)-\eta_{n}(y)|^{2}\abs{u_{n}(y)}^2}
{\abs{x-y}^{N+2s}}\dif x\dif y\Big)^{1/2}+\int_{\mathcal{B}_{3r_{n}}\setminus \mathcal{B}_{2r_{n}}}\abs{u_{n}}^{2}\dif x.
\end{aligned}\end{displaymath}
Arguing as the proof of \eqref{3.24}, we get
\begin{displaymath}\begin{aligned}
\int_{\mathbb{R}^{N}}\int_{\bbR^N}\frac{\abs{\eta_{n}(x)-\eta_{n}(y)}^{2}\abs{u_{n}(y)}^2}
{\abs{x-y}^{N+2s}}\dif x\dif y=o(1).
\end{aligned}\end{displaymath}
We thus deduce by \eqref{3.22} that
\begin{equation}\label{3.31}\begin{aligned}
\int_{\mathbb{R}^{N}}(-\Delta)^{s/2}v_{n}(-\Delta)^{s/2}w_{n}\dif x+\int_{\mathbb{R}^{N}}v_{n}w_{n}\dif x=o(1).
\end{aligned}\end{equation}
We now turn to the integration of the convolution term $F_n$. By combining \eqref{3.21} and the Hardy-Littlewood-Sobolev inequality, we derive that
\begin{equation}\label{3.32}\begin{aligned}
&\int_{\mathcal{B}_{3r_{n}}\setminus\mathcal{B}_{2r_{n}}}\int_{\mathbb{R}^{N}}F_n(x,y)\dif x\dif y\\
\leq&3\int_{\mathcal{B}_{3r_{n}}\setminus\mathcal{B}_{2r_{n}}}\int_{\mathbb{R}^{N}}
\frac{A_{\alpha}\abs{u_{n}(x)}^{p}\abs{u_{n}(y)}^{p}}{\abs{x-y}^{N-\alpha}}\dif x\dif y\\
\leq&C\Big(\int_{\mathcal{B}_{3r_{n}}\setminus\mathcal{B}_{2r_{n}}}
\abs{u_{n}(y)}^{\frac{2Np}{N+\alpha}}\dif y\Big)^{\frac{N+\alpha}{2N}}
\Big(\int_{\mathbb{R}^{N}}
\abs{u_{n}(x)}^{\frac{2Np}{N+\alpha}}\dif x\Big)^{\frac{N+\alpha}{2N}}=o(1).
\end{aligned}\end{equation}
Note that
\begin{equation}\label{3.33}\begin{aligned}
\int_{\mathcal{B}_{2r_{n}}}\int_{\mathbb{R}^{N}\setminus\mathcal{B}_{3r_{n}}}F_n(x,y)\dif x\dif y
=&\int_{\mathcal{B}_{2r_{n}}}\int_{\mathbb{R}^{N}\setminus\mathcal{B}_{3r_{n}}}
\frac{A_{\alpha}\abs{u_{n}(x)}^{p}\abs{u_{n}(y)}^{p}}{\abs{x-y}^{N-\alpha}}\dif x\dif y\\
&\leq\frac{C\|u_{n}\|_{L^{p}(\mathbb{R}^{N})}^{2}}{r^{N-\alpha}_{n}}
\leq\frac{C\|u_{n}\|_{2,s}^{2p}}{r^{N-\alpha}_{n}}=o(1).
\end{aligned}\end{equation}
We then  conclude that
\begin{equation}\label{3.34}\begin{aligned}
&\int_{\mathbb{R}^{N}}\int_{\bbR^N}F_n(x,y)\dif x\dif y\\
=&\int_{\mathcal{B}_{2r_{n}}}\int_{\mathcal{B}_{3r_{n}}\setminus\mathcal{B}_{2r_{n}}}F(x,y)\dif x\dif y
+\int_{\mathcal{B}_{2r_{n}}}\int_{\mathbb{R}^{N}\setminus\mathcal{B}_{3r_{n}}}F_n(x,y)\dif x\dif y\\
&+\int_{\mathcal{B}_{3r_{n}}\setminus\mathcal{B}_{2r_{n}}}\int_{\mathbb{R}^{N}}F_n(x,y)\dif x\dif y
+\int_{\mathbb{R}^{N}\setminus\mathcal{B}_{3r_{n}}}\int_{\mathcal{B}_{3r_{n}}}
F_n(x,y)\dif x\dif y\\
\leq&2\int_{\mathcal{B}_{2r_{n}}}\int_{\mathbb{R}^{N}\setminus\mathcal{B}_{3r_{n}}}F_n(x,y)\dif x\dif y
+3\int_{\mathcal{B}_{3r_{n}}\setminus\mathcal{B}_{2r_{n}}}\int_{\mathbb{R}^{N}}F_n(x,y)\dif x\dif y
=o(1).
\end{aligned}\end{equation}
Therefore, the claim \eqref{3.26} follows by a combination of \eqref{3.28}-\eqref{3.31} and \eqref{3.34}.

We now show the claim \eqref{3.27}. We define $\mathcal{K}_n:\bbR^N\times \bbR^N\to \bbR$ such that
\begin{displaymath}\begin{aligned}
\mathcal{K}_n(x,y)=\frac{\abs{u_{n}(x)}^{p}\abs{u_{n}(y)}^{p-2}u_{n}(y)v_{n}(y)-\abs{v_{n}(x)}^{p}\abs{v_{n}(y)}^{p}}{\abs{x-y}^{N-\alpha}}.
\end{aligned}\end{displaymath}
Arguing as the integration of $F_n(x,y)$, by \eqref{3.32} and \eqref{3.33}, we have
\begin{equation}\label{3.35}\begin{aligned}
&\int_{\mathbb{R}^{N}}\int_{\bbR^N}\mathcal{K}_n(x,y)\dif x\dif y\\
=&\int_{\mathbb{R}^{N}\setminus\mathcal{B}_{2r_{n}}}\int_{\mathcal{B}_{2r_{n}}}\mathcal{K}_n(x,y)\dif x\dif y
+\int_{\mathbb{R}^{N}}\int_{\mathcal{B}_{3r_{n}}\setminus\mathcal{B}_{2r_{n}}}\mathcal{K}_n(x,y)\dif x\dif y\\
\leq&3\int_{\mathbb{R}^{N}}\int_{\mathcal{B}_{3r_{n}}\setminus\mathcal{B}_{2r_{n}}}
\frac{ \abs{u_{n}(x)}^{p}\abs{u_{n}(y)}^{p}}{\abs{x-y}^{N-\alpha}}\dif x\dif y\\
\;&+\int_{\mathcal{B}_{2r_{n}}}\int_{\mathcal{B}_{3r_{n}}\setminus\mathcal{B}_{2r_{n}}}
\frac{ \abs{u_{n}(x)}^{p}\abs{u_{n}(y)}^{p}}{\abs{x-y}^{N-\alpha}}\dif x\dif y
=o(1).
\end{aligned}\end{equation}
Since $I'(u_{n})\rightarrow0$ in $(H^{s}_{G}(\mathbb{R}^{N}))^{\ast}$, and $v_{n}\in H^{s}_{G}(\mathbb{R}^{N})$, we have $\langle I'(u_{n}),v_{n}\rangle=o(1)$.
Then by \eqref{3.31} and \eqref{3.35}, we deduce that
\begin{equation}\label{3.36}\begin{aligned}
\int_{\mathbb{R}^{N}}\abs{(-\Delta)^{s/2}v_{n}}^{2}+\abs{v_{n}}^{2}\dif x=\int_{\mathbb{R}^{N}}\int_{\bbR^N}
\frac{A_{\alpha}\abs{v_{n}(x)}^{p}\abs{v_{n}(y)}^{p}}{\abs{x-y}^{N-\alpha}}\dif x\dif y+o(1).
\end{aligned}\end{equation}
Note that for any $i\neq m$,
\begin{displaymath}\begin{aligned}
\int_{\mathbb{R}^{N}}\abs{(-\Delta)^{s/2}v^{i}_{n}}^2
+\abs{v^{i}_{n}}^2\dif x=\int_{\mathbb{R}^{N}}\abs{(-\Delta)^{s/2}v^{m}_{n}}^2+\abs{v^{m}_{n}}^2 \dif x,
\end{aligned}\end{displaymath}
and
\begin{displaymath}\begin{aligned}
\int_{\mathbb{R}^{N}}\int_{\bbR^N}\frac{A_{\alpha}\abs{v^{i}_{n}(x)}^{p}\abs{v^{i}_{n}(y)}^{p}}{\abs{x-y}^{N-\alpha}}\dif x\dif y
=\int_{\mathbb{R}^{N}}\int_{\bbR^N}\frac{A_{\alpha}\abs{v^{m}_{n}(x)}^{p}\abs{v^{m}_{n}(y)}^{p}}{\abs{x-y}^{N-\alpha}}\dif x\dif y.
\end{aligned}\end{displaymath}
By combining \eqref{3.29}, \eqref{3.30} and \eqref{3.36}, we then infer that
\begin{displaymath}\begin{aligned}
\langle I'(v^{i}_{n}),v^{i}_{n}\rangle=o(1),\;\forall 1\leq i\leq|\mathcal{O}_{y_{n}}|.
\end{aligned}\end{displaymath}
By a similar argument as in \eqref{3.35}, we obtain
\begin{displaymath}\begin{aligned}
\int_{\mathbb{R}^{N}}\int_{\bbR^N}
\frac{ \abs{u_{n}(x)}^{p}\abs{u_{n}(y)}^{p-2}u_{n}(y)w_{n}(y)-\abs{w_{n}(x)}^{p}\abs{w_{n}(y)}^{p}}{\abs{x-y}^{N-\alpha}}\dif x\dif y
=o(1).
\end{aligned}\end{displaymath}
By taking advantage of $\lim_{n\to\infty}\langle I'(u_{n}),w_{n}\rangle=0$ and  \eqref{3.31}, we then conclude
$
\langle I'(w_{n}),w_{n}\rangle=o(1)
$.

\textbf{Case 2}: up to a subsequence, for any $\ell\geq 1$,  $\lim_{n\to\infty}\dist(y_n,\partial^\ell \mathcal D)=+\infty.$

In such case, up to a subsequence, we may assume for $z_i\neq z_j\in Gy_n$, there holds $\abs{x-y}\geq 4r_n$ for any $x\in B_{3r_n}(z_i)$ and $y\in B_{3r_n}(z_j)$.
We take the test function that
$$\phi_n^{G}=\sum_{g\in G}\xi(\abs{x-gy_n}/r_n)u_n.$$
We conclude that $\phi_n^G\in H^s_G(\bbR^N)$ since for every $\bar{g}\in G$,
$$\bar{g}\diamond \phi_n^{G}=\sum_{g\in G}\xi(\abs{x-g\bar{g}y_n}/r_n)u_n(\bar{g}^{-1}x)=\psi(\bar{g})\phi_n^G.$$
By repeating the arguments as above, we can also conclude \eqref{3.22}. Similarly, we define for each $z^j\in G{y_n}$ that
$$
v_n^j(x)=\eta(\abs{x-z^j}/r_n)u_n(x) \; \text{ and }\; w_n=u_n-\sum v_n^j.
$$
Since $\supp v_n^j\subset B_{3r_n}(z^j)$ for each $j$, we deduce that $v_n^j\in H^1(\bbR^N)$. It can also be checked that $w_n\in H^s_G(\bbR^N)$.
Similarly, we have the decompositions that
\begin{equation*}
\mathcal{E}(u_n)=\sum I(v_n^j)+I(w_n)+o(1), \;  \langle I'(v^{j}_{n}),v^{j}_{n}\rangle=o(1), \text{ and } \, \langle I'(w_{n}),w_{n}\rangle=o(1).
 \end{equation*}
 It then follows that there exists $t_n^j\in(0,+\infty)$ such that $t_n^jv_n^j\in \mathcal{N}$ with $t_n^j$ satisfying $\lim_{n\to\infty}t_n^j\to 1$ for each $1\leq j\leq \abs{G}$.
We therefrom deduce  that
\begin{align*}
c_G=\lim_{n\to\infty}I(u_n)&=\sum  \lim_{n\to\infty} I(v_n^j)+\lim_{n\to\infty} I(w_n)\geq
\sum  \lim_{n\to\infty}I(t_n^jv_n^j)\geq \abs{G}c_0.
\end{align*}
This is also in contradiction to the proposition \ref{Proposition 3.1} that $c_G<\abs{G}c_0$, so that \eqref{3.20} holds.

Therefore, when the compactness case happens, we take $a_{n}=(1-\Prj_{k})x_{n}\in\{0\}\times\mathbb{R}^{N-k}$ and
$T=M_{1}+R_{0}$, we get \eqref{3.17}; alternatively, if the dichotomy holds, we obtain \eqref{3.17} by choosing $a_{n}=(1-\Prj_{k})y_{n}\in\{0\}\times\mathbb{R}^{N-k}$
and $T=M_{2}+R_{1}$.

\section{Nodal structures for the saddle solutions}
\label{section4}
\vskip .1in

\quad\, We shall use the extension method of Caffarelli--Silvestre \cite{CS} to show the sign property of the saddle solution $u_{G}$.

Let $X^{1,s}(\mathbb{R}_{+}^{N+1})$ with $s\in(0,1)$ be the closure of  $C_{0}^{\infty}(\mathbb{R}_{+}^{N+1})$ under the norm
\begin{displaymath}
\begin{aligned}
\|U\|^{2}=\int_{\mathbb{R}_{+}^{N+1}}y^{1-2s}\abs{\nabla U}^{2} \dif x\dif y+\int_{\mathbb{R}^{N}}\abs{U(x,0)}^2\dif x.
\end{aligned}\end{displaymath}
For any $V\in X^{1,s}(\mathbb{R}_{+}^{N+1})$, $\Tr(V)$ denote its trace on $\mathbb{R}^{N}\times\{y=0\}$. Then
\begin{equation}\label{4.1}
\begin{aligned}
\|(-\Delta)^{s/2}\Tr(V)\|^{2}_{L^{2}(\mathbb{R}^{N})}\leq k^{-1}_{s}\int_{\mathbb{R}_{+}^{N+1}}y^{1-2s}\abs{\nabla V}^{2}\dif x\dif y,
\end{aligned}\end{equation}
where $k_{s}=2^{1-2s}\Gamma(1-s)\Gamma^{-1}(s)$.
For a given $u\in H^{s}(\mathbb{R}^{N})$,
there exists a unique $U\in X^{1,s}(\mathbb{R}_{+}^{N+1})$ such that
\begin{displaymath}
\begin{aligned}\left\{\begin{aligned}
   &-\Div(y^{1-2s}\nabla U)=0, && (x,y)\in\mathbb{R}_{+}^{N+1} ,& \\
     & U(x,0)=u(x),&& x \in \mathbb{R}^N,\\
  \end{aligned}
\right.
\end{aligned}\end{displaymath}
and
\begin{displaymath}
\begin{aligned}
\lim\limits_{y\rightarrow0}y^{1-2s}\frac{\partial U}{\partial y}(x,y)=-k_{s}(-\Delta)^{s}u(x).
\end{aligned}\end{displaymath}
We usually call $U=h_{s}(u)$ the $s$-harmonic extension of $u$. And it is given by
\begin{equation}\label{4.2}
\begin{aligned}\mathcal{F}(U)(\xi,y)=\mathcal{F}(u)(\xi)\psi(\abs{\xi}y),
\end{aligned}\end{equation}
where $\psi$ minimizes the functional
\begin{displaymath}
\begin{aligned}
H(\psi)=\int_{y>0}(\abs{\psi(y)}^{2}+\abs{\psi'(y)}^{2})y^{1-2s}\dif y.
\end{aligned}\end{displaymath}
Moreover, it is showed that
\begin{equation}\label{4.3}
\begin{aligned}
\int_{\mathbb{R}_{+}^{N+1}}y^{1-2s}\abs{\nabla h_{s}(u)}^{2}\dif x\dif y=\sqrt{k_{s}}\|(-\Delta)^{s/2}u\|_{L^{2}(\mathbb{R}^{N})}.
\end{aligned}\end{equation}

In what follows, the constant $k_{s}$ will be omitted for convenience.  For any $g\in G$, we define the action $g$ on $U\in X^{1,s}(\mathbb{R}^{N+1}_{+})$ by
\begin{displaymath}\begin{aligned}
g\circ U(x,y)=U(g^{-1}x,y),\;g^{-1}x=\diag\{g,1_{N-k}\}x.
\end{aligned}\end{displaymath}
Recall the unique epimorphism $\phi: G\to \{\pm 1\}$ and we define the space
\begin{displaymath}\begin{aligned}
X_{G}^{1,s}(\mathbb{R}_{+}^{N+1})=\{U\in X^{1,s}(\mathbb{R}_{+}^{N+1}):\,g\circ U(x,y)=\phi(g)U(x,y), \, \forall g\in G\}.
\end{aligned}\end{displaymath}
Let $J:X^{1,s}(\mathbb{R}^{N+1}_{+})\rightarrow\mathbb{R}$ be the functional defined by
\begin{displaymath}
\begin{aligned}
J(U)&=\frac{1}{2}\int_{\mathbb{R}_{+}^{N+1}}y^{1-2s}\abs{\nabla U}^{2}\dif x\dif y+\frac{1}{2} \int_{\mathbb{R}^{N}}\abs{U(x,0)}^{2} \dif x\\
&\quad -\frac{1}{2p}\int_{\mathbb{R}_+}\int_{\mathbb R^N} \frac{A_{\alpha}\abs{U(y,0)}^{p}\abs{U(x,0)}^{p}}{\abs{x-y}^{N-\alpha}}\dif x\dif y.
\end{aligned}\end{displaymath}
We then define analogously that
\begin{displaymath}\begin{aligned}
C_{G}=\inf_{U\in\mathcal{N}_{G}'}J(U),
\end{aligned}\end{displaymath}
where the constraint is
\begin{displaymath}
\begin{aligned}
\mathcal{N}_{G}'=\mathcal{N}'\cap X_{G}^{1,s}(\mathbb{R}_{+}^{N+1}),\end{aligned}\end{displaymath}
with Nehari manifold being
\begin{displaymath}\begin{aligned}
\mathcal{N}'=\big\{U\in X^{1,s}(\mathbb{R}_{+}^{N+1})\setminus\{0\} :
\|U\|^{2}=\int_{\mathbb{R}^{N}}(K_{\alpha}\ast\abs{U(\cdot,0)}^{p})\abs{U(x,0)}^{p}\dif x\big\}.
\end{aligned}\end{displaymath}

\begin{proposition}\label{Proposition 3.3} Let $u_{G}\in H^{s}_{G}(\mathbb{R}^{N})$ be a solution such that
\begin{displaymath}
\begin{aligned}
I(u_{G})=c_{G},\quad I'(u_{G})=0.
\end{aligned}\end{displaymath}
Then $u_{G}$ has a constant sign on $\mathcal{D}$ and hence $u_G$ has exactly $|G|$ nodal domains.
\end{proposition}
{\bf{ Proof}}.  We first show that $c_{G}=C_{G}$. On one hand, for any $u\in \mathcal{N}_{G}$, from \eqref{4.3}, we see that $U=h_{s}(u)\in X^{1,s}(\mathbb{R}^{N+1}_{+})$. Moreover,
by \eqref{4.2}, we conclude that $U\in\mathcal{N}'_{G}$. In fact, up to a constant, we have
 \begin{align*}
U(gx,y)&=\int_{\mathbb{R}^{N}}e^{2\pi \sqrt{-1} gx\cdot\xi}\int_{\mathbb{R}^{N}}e^{-2\pi \sqrt{-1} \xi\cdot z}u(z)\dif z\psi(\abs{\xi}y)\dif\xi\\
&=\int_{\mathbb{R}^{N}}e^{2\pi \sqrt{-1} x\cdot g^{-1}\xi}\int_{\mathbb{R}^{N}}e^{-2\pi\sqrt{-1} \xi \cdot  z} u(z)\dif z\psi(\abs{\xi}y)\dif\xi\\
&=\int_{\mathbb{R}^{N}}e^{2\pi \sqrt{-1} x\cdot\eta}\int_{\mathbb{R}^{N}}e^{-2\pi \sqrt{-1}g\eta  \cdot  z} u(z) \dif z\psi(\abs{g\eta}y)\dif \eta\\
&= \int_{\mathbb{R}^{N}}e^{2\pi \sqrt{-1} x\cdot\eta}\int_{\mathbb{R}^{N}}e^{-2\pi  \sqrt{-1}\eta  \cdot g^{-1}z} u(z) \dif z\psi(\abs{\eta}y)\dif \eta\\
&= \int_{\mathbb{R}^{N}}e^{2\pi \sqrt{-1} x\cdot\eta}\int_{\mathbb{R}^{N}}e^{-2\pi \sqrt{-1} \eta  \cdot \tilde{z}} u(g\tilde{z}) \dif \tilde{z}\psi(\abs{\eta}y)\dif \eta\\
&= \phi(g)\int_{\mathbb{R}^{N}}e^{2\pi \sqrt{-1} x\cdot\eta}\int_{\mathbb{R}^{N}}e^{-2\pi \sqrt{-1} \eta  \cdot \tilde{z}} u(\tilde{z}) \dif \tilde{z}\psi(\abs{\eta}y)\dif \eta
=\phi(g)U(x,y).
\end{align*}
We therefore deduce by the arbitrariness of $u\in \cN_G$ that
\begin{displaymath}\begin{aligned}
C_{G}\leq\inf_{u\in \mathcal{N}_{G}}J(h_{s}(u))=\inf_{u\in \mathcal{N}_{G}}I(u)=c_{G}.
\end{aligned}\end{displaymath}
On the other hand, for any $U\in \mathcal{N}'_{G}$,  there exists a unique $t_{u}>0$ such that $t_{u}U(x,0)\in\mathcal{N}_{G}$. We then have by \eqref{4.1} that

\begin{align*}
\inf\limits_{\mathcal{N}_{G}}I\leq & I(t_{u}U(x,0))
=\Big(\frac{1}{2}-\frac{1}{2p}\Big)\frac{ \Big(\displaystyle \int_{\mathbb{R}^{N}}\abs{(-\Delta)^{s/2}U(x,0)}^{2}
+\abs{U(x,0)}^{2}\dif x\Big)^{\frac{p}{p-1}}}{\Big(\displaystyle\int_{\mathbb{R}^{N}}(K_{\alpha}\ast\abs{U(\cdot,0)}^{p})
\abs{U(x,0)}^{p}\dif x\Big)^{\frac{1}{p-1}}}\\
\leq&\Big(\frac{1}{2}-\frac{1}{2p}\Big)\frac{\Big(\displaystyle\int_{\mathbb{R}^{N+1}_{+}}y^{1-2s}\abs{\nabla U}^{2}\dif x\dif y
+\int_{\mathbb{R}^{N}}\abs{U(x,0)}^{2}\dif x\Big)^{\frac{p}{p-1}}}{\Big(\displaystyle\int_{\mathbb{R}^{N}}(K_{\alpha}\ast\abs{U(\cdot,0)}^{p})
\abs{U(x,0)}^{p}\dif x\Big)^{\frac{1}{p-1}}}.
\end{align*}

Hence $c_{G}=\inf\limits_{\mathcal{N}_{G}}I\leq\inf\limits_{\mathcal{N}_{G}'}J=C_{G}.$

Let $U_{G}$ be the $s$-harmonic extension of $u_{G}\in H^{s}(\mathbb{R}^{N})$. Then
\begin{displaymath}\begin{aligned}
J(U_{G})=C_{G}=c_{G},\;\langle J'(U_{G}),\varphi\rangle =0, \;\forall\varphi\in X^{1,s}(\mathbb{R}^{N+1}_{+}).
\end{aligned}\end{displaymath}
Since $u_{G}\in C^{\beta}(\mathbb{R}^{N})$(see Lemma \ref{Lemma 3.1}), we have $U_{G}\in L^{\infty}(\overline{\mathbb{R}^{N+1}_{+}})$(see e.g., \cite[Lemma 4.4.]{CS1}) .
We now define $\overline{U}_G:\bbR_+^{N+1}\to \bbR$ such that
\begin{displaymath}
\begin{aligned}
\overline{U}_{G}=\left\{\begin{aligned}
 &\quad\abs{U_{G}(x,y)}, && (x,y)\in\overline{\mathcal{D}\times\{y\geq0\}} ,& \\
     &\phi(g)\abs{U_{G}(g^{-1}x,y)},&& (x,y)\in \overline{g(\mathcal{D})\times\{y\geq0\}}.\\
\end{aligned}
\right.
\end{aligned}\end{displaymath}
We then verify that $\overline{U}_{G}\in X^{1,s}_{G}(\mathbb{R}_{+}^{N+1})$.
Moreover, direct computations show us that
\begin{displaymath}
\begin{aligned}
J(\overline{U}_{G})=J(U_{G})=C_{G},\quad \langle J'(\overline{U}_{G}),\overline{U}_{G}\rangle =\langle J'(U_{G}),U_{G}\rangle=0.
\end{aligned}\end{displaymath}
Then similar arguments as that in \cite[Theorem 4.3]{W}, we conclude that $\overline{U}_{G}$ is a weak solution of

\begin{align}\label{eq4.4}
\left\{
  \begin{aligned}
     &-\Div(y^{1-2s}\nabla V)=0, && in\;\mathbb{R}_{+}^{N+1} ,& \\
     & \partial^{s}_{\nu}V=-V(x,0)
     +(K_{\alpha}\ast\abs{V(\cdot,0)}^{p})\abs{V(x,0)}^{p-2}V(x,0),
     && in\; \partial\mathbb{R}_{+}^{N+1},\\
  \end{aligned}
\right.\end{align}
where
 \begin{displaymath}
\begin{aligned}
\partial^{s}_{\nu}V=-\frac{1}{k_{s}}\lim\limits_{y\rightarrow0}y^{1-2s}
\frac{\partial V} {\partial y}(x,y)=(-\Delta)^{s} V(x,0).
\end{aligned}\end{displaymath}
Furthermore, we see that $\overline{u}_{G}=\overline{U}_{G}(x,0)$ is a weak solution of \eqref{eqfrc}.
By Lemma \ref{Lemma 3.1}, we have $\overline{u}_{G}\in C^{\beta}(\mathbb{R}^{N})$, then $(K_{\alpha}\ast\abs{\overline{u}_{G}}^{p})\abs{\overline{u}_{G}}^{p-2}
\overline{u}_{G}\in C^{\beta}(\mathbb{R}^{N})$. Since $U_{G}\in L^{\infty}(\mathbb{R}^{N+1}_{+})$, we have $\overline{U}_{G}\in L^{\infty}(\mathbb{R}^{N+1}_{+})$. Then applying a strong maximum principle (see e.g., \cite[Corollary 4.12]{CS1}) to the above system \eqref{eq4.4} in $\bbR_{+}^{N+1} \cap \mathcal{D}\times\{y\geq 0\}$,
we conclude  that $\overline{U}_{G}>0$ in $\mathcal{D}\times\{y\geq0\}$ and thus $\abs{u_{G}}=\overline{u}_{G}>0$ in $\mathcal{D}$.

\vskip .1in
\noindent{\bf Acknowledgement.} The authors would like to thank Professor Zhi-Qiang Wang for his useful discussions and encourages.
This paper is partially supported by NSFC (Grants 11771324, 11901456) and the Natural Science Foundation of Shaanxi Province
(Grant 2020JQ-120).


{\footnotesize

\begin {thebibliography}{99}


\bibitem{BBMP}
P. Belchior, H. Bueno, O. H. Miyagaki, G. A.  Pereira,
Remarks about a fractional Choquard equation: ground state, regularity and polynomial decay,
 \emph{Nonlinear Anal.}, \textbf{164} (2017)  38--53.

\bibitem{CS1}%
X. Cabr\'e, Y. Sire, Nonlinear equations for fractional Laplacians, I: Regularity, maximum principles, and Hamiltonian estimates,  \emph{Ann. Inst. H. Poincar\'e Anal. Non Lin\'aire}, \textbf{31}  (2014) 23--53.

\bibitem{Caf}%
L. Caffarelli, Nonlocal equations, drifts and games,  \emph{Nonlinear Partial Differential Equations, Abel
Symp.}, \textbf{7} (2012) 37--52.

\bibitem{CSS}%
L. Caffarelli, S. Salsa, L. Silvestre, Regularity estimates for the solution and the free boundary of the obstacle problem for the fractional Laplacian,  \emph{Invent. Math.}, \textbf{171}  (2008) 425--461.

\bibitem{CS}%
L. Caffarelli, L. Silvestre, An extension problem related to the fractional Laplacian,  \emph{Comm. Partial Differential Equations}, \textbf{32} (2007) 1245--1260.

\bibitem{CT2019}
S. Chen, X. Tang, Ground state solutions for asymptotically periodic fractional Choquard equations, \emph{
Electron. J. Qual. Theory Differ. Equ.}, (2019), Art. 2, 13 pp.
\bibitem{CL}%
Y. Chen, C. Liu, Ground state solutions for non-autonomous fractional Choquard equations,  \emph{Nonlinearity}, \textbf{29} (2016) 1827--1842.

\bibitem{CT}%
R. Cont, P. Tankov,  \emph{Financial Modelling with Jump Processes}, Chapman Hall/CRC
Financ. Math. Ser. Chapman and Hall/CRC, Boca Raton, FL (2004).

\bibitem{Cui2020}%
Y. -X. Cui, On nodal solutions of the fractional Choquard equation,  \emph{J. Math. Anal. Appl.}, \textbf{500} (2021), Art. 125152, 36pp.

\bibitem{DSS}%
P. d'Avenia, G. Siciliano, M. Squassina, On fractional Choquard equations,  \emph{Math. Models Methods Appl. Sci.}, \textbf{25}  (2015) 1447--1476.

\bibitem{Da}%
M. W. Davis,  \emph{The geometry and topology of Coxeter groups}, Vol 32, Princeton University Press, Princeton,
NJ., 2008.


\bibitem{DPV}%
 E. Di Nezza, G. Palatucci, E. Valdinoci, Hitchhiker's guide to the fractional Sobolev spaces,  \emph{Bull.
Sci. Math.}, \textbf{136} (2012) 521--573.


\bibitem{Le2009}
R. Frank, E. Lenzmann, On ground states for the $L^2$-critical boson star equation,
available at: \href{https://arxiv.org/abs/0910.2721}.

\bibitem{GMS}%
M. Ghimenti, V. Moroz, J.Van Schaftingen, Least action nodal solutions for the quadratic Choquard equation,  \emph{Proc. Amer. Math. Soc.}, \textbf{145}  (2017) 737--747.

\bibitem{GS}%
M. Ghimenti, J.Van Schaftingen, Nodal solutions for the Choquard equation,  \emph{J. Funct. Anal.}, \textbf{271}  (2016) 107--135.

\bibitem{GG}%
C. Gui, H. Guo, On nodal solutions of the nonlinear Choquard equation,  \emph{Adv. Nonlinear Stud.}, \textbf{19} (2019) 677--691.
\bibitem{GuoHu}
L. Guo, T. Hu, Existence and asymptotic behavior of the least energy solutions for fractional Choquard equations with potential well,  \emph{Math. Methods Appl. Sci.}, \textbf{41} (2018), 1145--1161.
\bibitem{HYY}
Z. Huang, J. Yang, W. Yu, Multiple nodal solutions of nonlinear Choquard equations,  \emph{Electron. J. Differential Equations}, (2017), Art. 268, 18 pp.

\bibitem{L}%
N. Laskin, Fractional quantum mechanics and L\'evy path integrals,  \emph{Phys. Lett. A}, \textbf{268}
(2000) 298--305.


\bibitem{Lieb}
E. H. Lieb,
   Existence and uniqueness of the minimizing solution of Choquard's nonlinear equation,
    \emph{Studies in Appl. Math.}, \textbf{57} (1976/77), 93--105.

\bibitem{LL}%
E. H. Lieb, M. Loss, \emph{ Analysis}, Graduate studies in mathematics, vol. 14. American Mathematical Society,
Providence (1997).

\bibitem{Li}%
P. -L. Lions, The Choquard equation and related questions,  \emph{Nonlinear Anal.}, \textbf{4}  (1980) 1063--1072.

\bibitem{Li2}%
P. -L. Lions, The concentration-compactness principle in the calculus of variations. The locally compact
case. I,   \emph{Ann. Inst. H. Poincar\'e Anal. Non Lin\'aire}, \textbf{1} (1984) 109--145.

\bibitem{LKZ}%
Q. Li, K. Teng, J. Zhang, Ground state solutions for fractional Choquard equations involving upper critical exponent,  \emph{Nonlinear Anal.}, \textbf{197} (2020), Art. 111846, 11 pp.

\bibitem{Luo}%
H. Luo, Ground state solutions of Poho\v{z}aev type for fractional Choquard equations with general nonlinearities,   \emph{Comput. Math. Appl.}, \textbf{77} (2019) 877--887.

\bibitem{MZ}%
L. Ma, L. Zhao, Classification of positive solitary solutions of the nonlinear Choquard equation,  \emph{Arch. Ration. Mech. Anal.}, \textbf{195} (2010) 455--467.

\bibitem{MZ2017}
P. Ma, J. Zhang,
Existence and multiplicity of solutions for fractional Choquard equations,  \emph{Nonlinear Anal.}, \textbf{164} (2017)  100--117.
\bibitem{MK2000}
R. Metzler, J. Klafter, The random walks guide to anomalous diffusion: A fractional dynamics approach,  \emph{Phys. Rep.}, \textbf{339} (2000) 1--77.

\bibitem{MoS}
V. Moroz, J. Van Schaftingen, Groundstates of nonlinear Choquard equations: existence, qualitative properties and decay asymptotics,  \emph{J. Funct. Anal.}, \textbf{265} (2013) 153--184.

\bibitem{MoS2}%
V. Moroz, J. Van Schaftingen, Existence of groundstates for a class of nonlinear Choquard equations,  \emph{Trans. Amer. Math. Soc.}, \textbf{367} (2015) 6557--6579.



\bibitem{MoS4}%
V. Moroz, J. Van Schaftingen, A guide to the Choquard equation,  \emph{J. Fixed Point Theory Appl.}, \textbf{19}  (2017) 773--813.

\bibitem{P}%
S. Pekar,  \emph{Untersuchung uber die Elektronentheorie der Kristalle}, Akademie Verlag, Berlin, 1954.

\bibitem{Pe}%
R. Penrose, Quantum computation, entanglement and state reduction,  \emph{R. Soc. Lond. Philos. Trans. Ser. A Math. Phys. Eng. Sci.}, \textbf{356} (1998) 1927--1939.

\bibitem{SGY}%
Z. Shen, F. Gao, M. Yang, Ground states for nonlinear fractional Choquard equations with general nonlinearities,  \emph{Math. Methods Appl. Sci.}, \textbf{39}  (2016) 4082--4098.

\bibitem{SX1}%
J. Van Schaftingen, J. Xia, Choquard equations under confining external potentials, \emph{ NoDEA Nonlinear Differential Equations Appl.}, \textbf{24} (2017), Art. 1, 24 pp.


\bibitem{Th}%
A. Thomas,  \emph{Geometric and topological aspects of Coxeter groups and buildings}, Vorlesungsmanuskript, 2017.


\bibitem{XW2}%
Z.-Q. Wang, J. Xia, Saddle solutions for the Choquard equation \uppercase\expandafter{\romannumeral2},
 \emph{Nonlinear, Anal.}, \textbf{201} (2020), Art. 112053, 25pp.

\bibitem{Weth2006} T. Weth,
     Energy bounds for entire nodal solutions of autonomous superlinear equations,
       \emph{Calc. Var. Partial Differential Equations}, \textbf{27} (2006) 421--437.

\bibitem{W}%
M. Willem, \emph{ Minimax Theorems}, Progress in nonlinear differential equations and their Applications, vol.24, Birkh$\acute{a}$user, Boston, Mass, 1996.

\bibitem{XW1}%
J. Xia, Z.-Q. Wang, Saddle solutions for the Choquard equation,
 \emph{Calc. Var. Partial Differential Equations}, \textbf{58}  (2019), Art. 85, 30 pp.

\bibitem{XiaZhang2020}%
J. Xia, X. Zhang, Saddle solutions for the critical Choquard equation,   \emph{ Calc. Var. Partial Differential Equations}, \textbf{60} (2021), Art. 53, 29 pp.

\bibitem{Xia2021}%
J. Xia, Saddle solutions for the Choquard equation with a general nonlinearity, (2020), submitted.

\bibitem{ZW}%
W. Zhang, X. Wu, Nodal solutions for a fractional Choquard equation,  \emph{J. Math. Anal. Appl.}, \textbf{464}  (2018) 1167--1183.




\end {thebibliography}
}
\end{document}